\documentclass[a4paper]{article}
\usepackage{amssymb,amscd,amsfonts,mathrsfs,amsmath}
\usepackage[all]{xy}
\input amssym.def

\def\double{\mathbb}

\def\cc{{\double C}}
\def\nn{{\double N}}
\def\zz{{\double Z}}

\newtheorem{theorem}{Theorem}[section]
\newtheorem{lemma}[theorem]{Lemma}
\newtheorem{corollary}[theorem]{Corollary}
\newtheorem{definition}[theorem]{Definition}
\newtheorem{proposition}[theorem]{Proposition}
\newtheorem{remark}[theorem]{Remark}

\def\res{\mathop{\mathrm{Res}}\limits_{z=0}}
\def\Pf{\mathop{\mathrm{Pf}}}

\def\si{\sigma}

\def\cinf{C^{\infty}}
\def\cinfc{C^{\infty}_c}

\newcommand{\be}{\begin{equation}}
\newcommand{\ee}{\end{equation}}
\newcommand{\beq}{\begin{eqnarray}}
\newcommand{\eeq}{\end{eqnarray}}
\newcommand{\om}{\omega}
\newcommand{\Om}{\Omega}
\newcommand{\al}{\alpha}
\def\nat{\natural}

\newcommand{\La}{\Lambda}
\newcommand{\la}{\lambda}
\newcommand{\Ec}{{\mathscr E}}

\newcommand{\non}{\nonumber}
\newcommand{\eps}{\varepsilon}

\newcommand{\Rc}{{\mathscr R}}
\newcommand{\Mc}{{\mathscr M}}
\newcommand{\Nc}{{\mathscr N}}
\newcommand{\Ic}{{\mathscr I}}
\newcommand{\Jc}{{\mathscr J}}

\newcommand{\Ind}{{\mathop{\mathrm{Ind}}}}
\def\ch{\mathrm{ch}}

\def\cs{\mathrm{cs}}
\def\re{\mathrm{Re}}

\newcommand{\Tr}{{\mathop{\mathrm{Tr}}}}

\newcommand{\tr}{{\mathop{\mathrm{tr}}}}
\newcommand{\Ac}{{\mathscr A}}
\newcommand{\Gc}{{\mathscr G}}
\newcommand{\te}{\theta}

\newcommand{\cqfd}{\hfill\rule{1ex}{1ex}}

\def\Id{\mathrm{Id}}

\def\d{\partial}
\def\dd{\mathrm{\bf d}}

\def\Hc{{\mathscr H}}

\def\Bc{{\mathscr B}}
\def\Cc{{\mathscr C}}
\def\Jc{{\mathscr J}}

\def\Fc{{\mathscr F}}
\def\Pc{{\mathscr P}}
\def\Qc{{\mathscr Q}}

\def\ker{\mathop{\mathrm{Ker}}}

\def\bb{\overline{b}}

\def\sib{\overline{\si}}

\def\Tt{\widetilde{T}}

\def\gt{\tilde{g}}

\def\pit{\tilde{\pi}}
\def\Omh{\widehat{\Omega}}
\def\Th{\widehat{T}}
\def\chih{\widehat{\chi}}
\def\etah{\widehat{\eta}}

\def\Rch{\widehat{\mathscr R}}
\def\Gch{\widehat{\mathscr G}}
\def\Mch{\widehat{\mathscr M}}
\def\Mct{\widetilde{\mathscr M}}

\def\Omc{\Omega_c}

\def\mod{\ \mathrm{mod}\ }

\def\eh{\hat{e}}
\def\zh{\hat{z}}
\def\gh{\hat{g}}

\def\bint{- \mspace{-21mu} \int}

\begin{document}

\begin{center}

{\bf EXTENSIONS AND RENORMALIZED TRACES}
\vskip 1cm
{\bf Denis PERROT}
\vskip 0.5cm
Universit\'e de Lyon, Universit\'e Lyon 1,\\
CNRS, UMR 5208 Institut Camille Jordan,\\
43, bd du 11 novembre 1918, 69622 Villeurbanne Cedex, France \\[2mm]
{\tt perrot@math.univ-lyon1.fr}\\[2mm]
\today
\end{center}
\vskip 0.5cm
\begin{abstract}
It has been shown by Nistor \cite{Ni} that given any extension of associative algebras over $\cc$, the connecting morphism in periodic cyclic homology is compatible, under the Chern-Connes character, with the index morphism in lower algebraic $K$-theory. The proof relies on the abstract properties of cyclic theory, essentially excision, which does not provide explicit formulas a priori. Avoiding the use of excision, we explain in this article how to get explicit formulas in a wide range of situations. The method is connected to the renormalization procedure introduced in our previous work on the bivariant Chern character for quasihomomorphisms \cite{P5, P6}, leading to ``local'' index formulas in the sense of non-commutative geometry. We illustrate these principles with the example of the classical family index theorem: we find that the characteristic numbers of the index bundle associated to a family of elliptic pseudodifferential operators are expressed in terms of the (fiberwise) Wodzicki residue. 
\end{abstract}

\vskip 0.5cm


\section{Introduction}

Some years ago Cuntz and Quillen were able to show that excision holds in complete generality for periodic cyclic (co)homology of associative algebras \cite{CQ97}. That is, given any extension (short exact sequence) of algebras over $\cc$,
\be
(E)\ :\ 0 \to \Bc \to \Ec \to \Ac \to 0
\ee
there exists an associated six-term exact sequence relating the periodic cyclic homology of $\Bc$, $\Ec$, $\Ac$, and similarly for cohomology. Using the abstract properties of the theory, Nistor \cite{Ni} then proved that the connecting morphism $HP_1(\Ac)\to HP_0(\Bc)$ of the cyclic homology exact sequence is compatible, via the Chern-Connes character, with the \emph{index map} induced by the extension $(E)$ on algebraic $K$-theory in low degrees \cite{M}:
\be
\Ind_E:\ K_1(\Ac)\to K_0(\Bc)\ . \label{ind}
\ee
In principle this allows to state a general ``higher index theorem'', in the sense that the pairing of any periodic cyclic cohomology class $[\tau]\in HP^0(\Bc)$ with the image of (\ref{ind}) can be computed as the pairing of its boundary $E^*([\tau])$ with $K_1(\Ac)$. Here $E^*:HP^0(\Bc)\to HP^1(\Ac)$ denotes the connecting morphism in cohomology. Nistor proves this theorem first in the case of a universal extension, for which the periodic cyclic cohomology is simply represented by traces over $\Bc$, then the general case follows from the naturality of the index morphism in $K$-theory and the naturality of the boundary map in periodic cyclic cohomology. Although very elegant and general, this proof does not provide explicit formulas for the cocycle $E^*([\tau])$. In principle the proof of excision in \cite{CQ97} should lead to explicit formulas, but they turn out to be extremely complicated in general, and moreover are not \emph{local} in contrast with, for instance, the residue index formula of Connes and Moscovici \cite{CM95}. \\

The goal of the present article is to present an explicit construction of the connecting morphism $E^*$ avoiding as much as possible the use of excision, and giving an alternative (direct) proof of Nistor's index theorem. One knows from the work of Cuntz and Quillen \cite{CQ95} that any cyclic cohomology class $[\tau]\in HP^0(\Bc)$ can be represented by a trace over an adequate extension $0\to\Jc\to\Rc\to\Bc\to 0$ of $\Bc$, or equivalently by a trace over some power of this extension (think for example about the operator trace on a Schatten ideal). Our basic observation is the following: if the extensions $0\to\Jc\to\Rc\to\Bc\to 0$ and $0\to\Bc\to\Ec\to\Ac\to 0$ fit together in a commutative diagram (see (\ref{diag})), then $E^*([\tau])$ is explicitely given by a fairly simple formula based on a ``renormalization'' procedure explained in section \ref{scm}. The proof that actually \emph{any} cyclic cohomology class over $\Bc$ can be represented in this way requires the knowledge of excision. Fortunately many cyclic cohomology classes appear naturally equipped with the required diagram, so we are able to circumvent excision completely in this situation. \\
Let us mention that the term ``renormalization'' is inspired by our previous work on the bivariant Chern character for quasihomomorphisms \cite{P5, P6, P7}, where it was argued that this procedure yields local index formulas automatically. This is related to the well-known anomalies of quantum field theory \cite{P4}. In fact we show in section \ref{squa} that when the extension $(E)$ is \emph{invertible} in a specific sense (Definition \ref{dinv}), the map $E^*$ coincides with the bivariant Chern character of the odd quasihomomorphism associated to the extension. This allows to give an alternative proof of Nistor's index theorem in section \ref{sind}: for any $[\tau]\in HP^0(\Bc)$ and $[g]\in K_1(\Ac)$, one has the equality of pairings
\be
\langle [\tau], \Ind_E([g]) \rangle = \sqrt{2\pi i}\, \langle E^*([\tau]), [g] \rangle \label{the}
\ee
that is, the index map is adjoint to the connecting morphism in periodic cyclic cohomology. The overall factor $\sqrt{2\pi i}$ comes from our particular choice of normalization for the pairings between cyclic cohomology and $K$-theory: just note that this choice is the only one compatible with the bivariant Chern character and Bott periodicity for topological algebras. Since we will consider only algebras without additional structure in this article, this factor is irrelevant. The index theorem is shown in two steps: first we reduce to the case of an invertible extension, and then (\ref{the}) is the consequence of an explicit computation. Thus in contrast with \cite{Ni}, excision is not directly used in the proof. \\
In section \ref{spseu} we show on the example of the family index theorem that our construction of $E^*([\tau])$ effectively leads to local formulas. Thus we consider a proper submersion of smooth manifolds without boundary $M\to B$. A canonical extension $(E)$ is obtained by taking $\Ec$ as the algebra of smooth families of (fiberwise) classical pseudodifferential operators of order zero, $\Bc$ as the ideal of order $-1$ pseudodifferential operators, and $\Ac$ as the commutative algebra of smooth functions over the cotangent sphere bundle of the fibers. The projection $\Ec\to\Ac$ thus carries a family of pseudodifferential operators to its family of leading symbols. Then any de Rham cycle in the base manifold $B$ gives rise to a cyclic cocycle $\tau$ over the algebra $\Bc$; notice however that this requires to choose a connection on the submersion. Using zeta-function renormalization, we find that the cyclic cohomology class $E^*([\tau])\in HP^0(\Ac)$ is given explicitely in terms of a fiberwise Wodzicki residue applied to some families of pseudodifferential operators, involving the connection and its curvature. Interestingly, the formula is a higher analogue of the famous Radul cocycle \cite{Ra}. Finally if $Q$ is a family of elliptic pseudodifferential operators with symbol class $[g]\in K_1(\Ac)$, the pairing between $[\tau]$ and the ``index bundle'' $\Ind_E([g])\in K_0(\Bc)$ is the evaluation of this higher Radul cocycle on certain polynomials built from $Q$, its parametrix $P$, and the connection.

\section{Connecting morphism}\label{scm}

Let us recall the Cuntz-Quillen formalism of cyclic cohomology \cite{CQ95}, since it is particularly well-adapted to extensions. The basic fact is that any cyclic cohomology class of even degree over an associative algebra $\Bc$ can be represented by a trace over some extension $\Rc$
\be
0\to\Jc\to\Rc\to\Bc\to 0
\ee
vanishing on the large powers of the ideal $\Jc$. A cyclic cohomology class of odd degree over $\Bc$ can be represented by a cyclic one-cocycle on $\Rc$ with similar vanishing properties. This motivates the definition of the $X$-complex of any algebra $\Rc$. It is the $\zz_2$-graded complex
\be
X(\Rc):\quad \Rc\ \mathop{\rightleftarrows}^{\nat \dd}_{\bb}\ \Om^1\Rc_{\nat}\ , 
\ee
where $\Om^1\Rc_{\nat}=\Om^1\Rc/[\Rc,\Om^1\Rc]$ is the quotient of the $\Rc$-bimodule of universal one-forms by its commutator subspace. The class of a generic element $(x_0\dd x_1 \mod [,])\in \Om^1\Rc_{\nat}$ is usually denoted by $\nat x_0\dd x_1$. The map $\nat \dd:\Rc\to \Om^1\Rc_{\nat}$ thus sends $x\in\Rc$ to $\nat \dd x$. Also, the Hochschild boundary map $b:\Om^1\Rc\to \Rc$ vanishes on the commutator subspace, hence passes to a well-defined map $\bb:\Om^1\Rc_{\nat}\to\Rc$. Explicitly the image of $\nat x_0\dd x_1$ by $\bb$ is the commutator $[x_0,x_1]$. These maps satisfy $\nat \dd\circ \bb=0$ and $\bb\circ\nat \dd=0$, so that $X(\Rc)$ endowed with the boundary operator $\d = \nat\dd\oplus\bb$ indeed defines a $\zz_2$-graded complex.\\
If $\Jc\subset \Rc$ is a two-sided ideal, Cuntz and Quillen define a decreasing filtration of $X(\Rc)$ by the following subcomplexes indexed by integers $n\in\zz$
\beq 
F_{\Jc}^{2n}X(\Rc) &:& \Jc^{n+1}+[\Jc^n,\Rc] \ \rightleftarrows \ \nat \Jc_{(+)}^n\dd\, \Rc \label{filtration}\\
F_{\Jc}^{2n+1}X(\Rc) &:& \Jc^{n+1}\ \rightleftarrows\ \nat(\Jc_{(+)}^{n+1}\dd\, \Rc + \Jc_{(+)}^n\dd\, \Jc)\ ,\non
\eeq
where $\Jc_{(+)}^n$ is equal to the power $\Jc^n$ for $n>0$ and equal to the unitalized algebra $\Rc^+=\Rc\oplus\cc$ for $n\leq 0$. The $\Jc$-adic completions of the algebra $\Rc$ and of the complex $X(\Rc)$ are defined as projective limits
\be
\Rch=\varprojlim_n \Rc/\Jc^n\ ,\qquad X(\Rch)=\varprojlim_n X(\Rc)/F^n_{\Jc}X(\Rc)
\ee
where $\Rch$ is viewed as a \emph{pro-algebra} and $X(\Rch)$ as a \emph{pro-complex} \cite{CQ97}. It follows that any cocycle $\tau:X(\Rch)\to\cc$ represents a cyclic cohomology class over $\Bc$ (here $\tau$ is viewed as a linear map between pro-complexes, that is, a linear map on $X(\Rc)$ vanishing on $F^n_{\Jc}X(\Rc)$ for some $n\gg 0$). In particular let $T\Bc=\Bc\oplus \Bc\otimes\Bc\oplus\ldots$ be the non-unital tensor algebra and denote by $J\Bc$ the kernel of the multiplication homomorphism $T\Bc\to\Bc$. Then the extension $0\to J\Bc\to T\Bc\to\Bc\to 0$ is universal among all extensions $0\to\Jc\to\Rc\to\Bc\to 0$ in the sense that one has a classifying homomorphism $T\Bc\to\Rc$ defined up to homotopy, which restricts to a homomorphism $J\Bc\to \Jc$. Thus any cocycle over $X(\Rch)$ can be pulled back to a cocycle over $X(\Th\Bc)$. In particular the cohomology group $H^*(X(\Th\Bc))$ is isomorphic to the periodic cyclic cohomology $HP^*(\Bc)$, see \cite{CQ95}.\\

It follows from the proof of excision in periodic cyclic cohomology \cite{CQ97}, that any extension
\be
(E)\ :\ 0\to\Bc\to\Ec\to\Ac\to 0
\ee
gives rise to a connecting morphism $HP^i(\Bc)\to HP^{i+1}(\Ac)$, $i\in \zz_2$. Here we shall present a way to calculate the connecting morphism, assuming that the cyclic cohomology classes of $\Bc$ are put into a suitable form. They will be represented not only by traces over some extension $\Rc$, but more generally by traces over some power $\Rc^n$: indeed cyclic cohomology classes often arise as traces over finitely summable operator ideals \cite{C86}. Since one has to choose extensions of both algebras $\Bc$ and $\Ac$ to represent their cyclic cohomology, we first define the notion of a lifting for the extension $(E)$:

\begin{definition}\label{dlift}
We say that an extension $0\to\Rc\to\Mc\to\Pc\to 0$ is a lifting of $0\to\Bc\to\Ec\to\Ac\to 0$, if both fit into a commutative diagram
\be
\vcenter{\xymatrix{
 & 0 \ar[d] & 0 \ar[d] & 0 \ar[d] &  \\
0 \ar [r] & \Jc \ar[r] \ar[d] & \Nc \ar[r] \ar[d] & \Qc  \ar[r] \ar[d] & 0  \\
0 \ar [r] & \Rc \ar[r] \ar[d] & \Mc \ar[r] \ar[d] & \Pc  \ar[r] \ar[d] & 0  \\
0 \ar [r] & \Bc \ar[r] \ar[d] & \Ec \ar[r] \ar[d] & \Ac  \ar[r] \ar[d] & 0  \\
 & 0 & 0 & 0 &  }} \label{diag}
\ee
where all rows and columns are extensions.
\end{definition}

Morally the columns of (\ref{diag}) will be used to represent cyclic cohomology classes of $\Bc$, $\Ec$, $\Ac$ respectively. Because the central algebra $\Mc$ has two distinguished ideals $\Rc$ and $\Nc$, there are several ways to filter the complex $X(\Mc)$. First we focus on the middle column, i.e. the extension $0\to \Nc\to \Mc\to\Ec\to 0$. We denote the $\Nc$-adic completions with a hat:
\be
\Mch=\varprojlim_n \Mc/\Nc^n\ ,\qquad X(\Mch)=\varprojlim_n X(\Mc)/F^n_{\Nc}X(\Mc)\ .
\ee
There is a second extension $0\to\Rc+ \Nc\to \Mc\to\Ac\to 0$ associated to the diagonal of (\ref{diag}). The corresponding $(\Rc+\Nc)$-adic completions will be denoted with a tilde. Of course this completion of $X(\Mc)$ could be defined via the filtration by the subcomplexes $F^n_{\Rc+\Nc}X(\Mc)$, but there is an equivalent construction starting from the above pro-complex $X(\Mch)$ filtered by the subcomplexes $F^n_{\Rc}X(\Mch)= \varprojlim_k F^n_{\Rc}X(\Mc)/(F^n_{\Rc}X(\Mc)\cap F^k_{\Nc}X(\Mc))$:
\be
\Mct=\varprojlim_n \Mc/(\Rc+\Nc)^n\ ,\qquad X(\Mct)=\varprojlim_n X(\Mch)/F^n_{\Rc}X(\Mch)\ .
\ee
The cocycles over $F^n_{\Rc}X(\Mch)$ are chain maps $\tau:F^n_{\Rc}X(\Mch)\to\cc$. We will show below that they represent cyclic cohomology classes over $\Bc$. Let us describe now how to compute the connecting morphism associated to the initial extension $(E):0\to\Bc\to\Ec\to\Ac\to 0$. 

\begin{definition}\label{dren}
Let $\tau$ be a cocycle over the subcomplex $F^n_{\Rc}X(\Mch)$ for some $n\geq 1$. A \emph{renormalization} of $\tau$ is a linear map $\tau_R:X(\Mch)\to\cc$ extending $\tau$.
\end{definition}
Of course $\tau_R$ is usually not a cocycle over $X(\Mch)$. Its coboundary $\tau_R\d$ is however a cocycle vanishing on $F^n_{\Rc}X(\Mch)$ by construction. Here $\d=\nat\dd\oplus\bb$ is the $X$-complex boundary map. Since $F^{n+k}_{\Rc}X(\Mch)\subset F^n_{\Rc}X(\Mch)$ for any $k\geq 0$, one sees that $\tau_R\d$ descends to a unique cocycle over $X(\Mct)$. It remains to pull it back to $X(\Th\Ac)$ in order to get a periodic cyclic cohomology class over $\Ac$. Choose a linear splitting $\si:\Ac\to \Mc$ of the diagonal homomorphism $\Mc\to \Ac$. The universal property of the tensor algebra $T\Ac$ allows to extend $\si$ to a homomorphism $\si_*:T\Ac\to \Mc$ by setting $\si_*(a_1\otimes\ldots\otimes a_k)=\si(a_1)\ldots\si(a_k)$. Then $\si_*$ sends the ideal $J\Ac=\ker(T\Ac\to\Ac)$ to the ideal $\Rc+\Nc=\ker(\Mc\to\Ac)$. This may be depicted through the following commutative diagram where all arrows except the dashed one are homomorphisms of algebras:
\be
\vcenter{\xymatrix{
0 \ar[r] & J\Ac \ar[r] \ar[d]_{\si_*} & T\Ac \ar[r] \ar[d]_{\si_*} & \Ac \ar[r] \ar@{=}[d] \ar@{.>}[dl]^{\si} & 0 \\
0\ar[r] & \Rc + \Nc \ar[r]  & \Mc \ar[r] & \Ac \ar[r]  & 0 }}
\ee
Consequently $\si_*$ extends to a homomorphism of pro-algebras $\Th\Ac\to \Mct$. This in turn induces a chain map $\si_*:X(\Th\Ac)\to X(\Mct)$.

\begin{proposition}\label{pren}
Consider an extension $(E):0\to\Bc\to\Ec\to\Ac\to 0$ with lifting $0\to\Rc\to\Mc\to\Pc\to 0$ as in diagram (\ref{diag}). The map sending a cocycle $\tau$ over $F^n_{\Rc}X(\Mch\,)$ to the cocycle $\tau_R\d\circ\si_*$ over $X(\Th\Ac)$, for an arbitrary choice of renormalization $\tau_R$, descends to a morphism in cohomology
\be
E^n:H^i\big(F^n_{\Rc}X(\Mch)\big) \to HP^{i+1}(\Ac)\ , \quad i\in \zz_2\ . \label{mor}
\ee
The latter does not depend on the choice of renormalization, nor on the linear splitting $\si:\Ac\to\Mc$. If the exact sequence $(E)$ is split by a homomorphism $\Ac\to\Ec$, then $E^n$ vanishes. Finally, $E^{n+1}$ composed with the natural pullback $H^*\big(F^n_{\Rc}X(\Mch)\big)\to H^*\big(F^{n+1}_{\Rc}X(\Mch)\big)$ coincides with $E^n$.
\end{proposition}
{\it Proof:} If $\tau_R'$ is another choice of linear extension for $\tau$, the difference $\tau_R'-\tau_R$ vanishes on $F^n_{\Rc}X(\Mch)$, hence descends to a cochain over $X(\Mct)$. The difference $\tau_R'\d-\tau_R\d=(\tau_R'-\tau_R)\d$ is therefore a coboundary over $X(\Mct)$. The cyclic cohomology class of $\tau_R\d\circ\si_*$ is thus renormalization-independent.\\
If $\tau=\varphi\d$ is the coboundary of a cochain $\varphi$ over $F^n_{\Rc}X(\Mch)$, then one can extend $\varphi$ to a cochain $\varphi_R$ over $X(\Mch)$ and take $\tau_R=\varphi_R\d$. Then $\tau_R\d=0$ and $E^n$ is well-defined in cohomology. \\
As observed by Cuntz and Quillen \cite{CQ95}, two different choices of linear splittings $\si:\Ac\to\Mc$ induce homotopic homomorphisms $\si_*:\Th\Ac\to\Mct$, and the resulting chain maps $X(\Th\Ac)\to X(\Mct)$ are homotopy equivalent. Thus the cyclic cohomology class of $\tau_R\d\circ\si_*$ is independent of $\si$.\\
Suppose that $(E)$ is split by a homomorphism $\rho:\Ac\to\Ec$. Then choose any linear splitting $\ell:\Ec\to\Mc$ of the projection homomorphism $\Mc\to\Ec$ and put $\si=\ell\circ\rho$. By the universality of $T\Ac$ one gets a commutative diagram
$$
\vcenter{\xymatrix{
0 \ar[r] & J\Ac \ar[r] \ar[d]_{\si_*} & T\Ac \ar[r] \ar[d]_{\si_*} & \Ac \ar[r] \ar[d]^{\rho} \ar@{.>}[dl]^{\si} & 0 \\
0\ar[r] & \Nc \ar[r]  & \Mc \ar[r] & \Ec \ar[r] \ar@{.>}@/^1pc/[l]^{\ell} & 0 }}
$$
Since $J\Ac$ lands in $\Nc$, the homomorphism $\si_*:\Th\Ac\to\Mct$ actually factors through $\Mch$. Hence the composite map $\tau_R\circ\si_*$ is a well-defined cochain over $X(\Th\Ac)$, the cocycle $\tau_R\d\circ\si_*=(\tau_R\circ\si_*)\d$ is a coboundary, and $E^n$ vanishes. The last assertion is obvious.  \cqfd\\

\begin{remark}\label{rem}\textup{If $\Nc$ is nilpotent then the projective limit $\Mch$ reduces to $\Mc$. This has the following important consequence concerning the cocycles of \emph{even} degree $\tau : F^{2n+1}_{\Rc}X(\Mc)\to\cc$. Indeed, such a cocycle is a linear map $\tau:\Rc^{n+1}\to\cc$ vanishing on the commutator subspace $[\Rc^{n+1},\Mc]+[\Rc^n,\Rc]$. In particular for $n\geq 1$ the inclusion $[\Rc^{n+1},\Mc]\subset [\Rc^n,\Rc]$ holds, hence \emph{any reference to $\Mc$ disappears}. We can conclude that a cohomology class in $H^0(F^{2n+1}_{\Rc}X(\Mc))$ is simply represented by a trace over the $(n+1)$-th power of the nilpotent extension $0\to\Jc\to\Rc\to\Bc\to 0$, that is, a linear map $\tau:\Rc^{n+1}\to\cc$ vanishing on $[\Rc^n,\Rc]$ for $n\gg 0$.\\
In the general case $\Nc$ is not nilpotent, and an even cocycle $\tau$ has to verify the additional condition that it vanishes on $\Nc^k$ for some $k\gg 0$. However one recovers the nilpotent situation after replacing the second row of (\ref{diag}) by the new extension $0\to\Rc/(\Rc\cap\Nc^k)\to\Mc/\Nc^k\to\Pc/\Qc^k\to 0$, and the first row by $0\to\Jc/(\Jc\cap\Nc^k)\to\Nc/\Nc^k\to\Qc/\Qc^k\to 0$. Then $\tau$ still defines a cocycle for this diagram and the new ideal $\Nc/\Nc^k$ is nilpotent. We nevertheless prefer to stay in the general context since important examples of universal extensions are not nilpotent.}
\end{remark}

Given any extension $(E):0\to\Bc\to\Ec\to\Ac\to 0$ there always exists a lifting in the sense of Definition \ref{dlift}. Indeed one can consider the following \emph{universal lifting}
\be
\vcenter{\xymatrix{
 & 0 \ar[d] & 0 \ar[d] & 0 \ar[d] &  \\
0 \ar [r] & J(\Bc:\Ec) \ar[r] \ar[d] & J\Ec \ar[r] \ar[d] & J\Ac  \ar[r] \ar[d] & 0  \\
0 \ar [r] & T(\Bc:\Ec) \ar[r] \ar[d] & T\Ec \ar[r] \ar[d] & T\Ac  \ar[r] \ar[d] & 0  \\
0 \ar [r] & \Bc \ar[r] \ar[d] & \Ec \ar[r] \ar[d] & \Ac  \ar[r] \ar[d] & 0  \\
 & 0 & 0 & 0 &  }} \label{univ}
\ee
where the ideal $T(\Bc:\Ec)$ (resp. $J(\Bc:\Ec)$) denotes the kernel of the homomorphism $T\Ec\to T\Ac$ (resp. $J\Ec\to J\Ac$). The universal property of the tensor algebras $T\Ec$ and $T\Ac$ induce classifying maps from the second and third column of (\ref{univ}) to the second and third column of (\ref{diag}) respectively, and this in turn implies a classifying map for the first column also. The central classifying homomorphism $T\Ec\to\Mc$ and all other ones are defined up to homotopy, which ensures a canonical pullback morphism $H^*\big(F^n_{\Rc}X(\Mch)\big)\to H^*\big(F^n_{T(\Bc:\Ec)}X(\Th\Ec)\big)$ for all $n$. In fact the excision property of periodic cyclic cohomology \cite{CQ97} shows the following

\begin{lemma} Let $(E): 0\to\Bc\to\Ec\to\Ac\to 0$ be an extension. Then for any $n\geq 1$ one has an isomorphism 
\be
H^*\big(F^n_{T(\Bc:\Ec)}X(\Th\Ec)\big)\cong HP^*(\Bc)\ .
\ee
\end{lemma}
{\it Proof:} According to the terminology of Cuntz and Quillen the pro-algebra $\Th\Ec$ is a quasi-free extension of $\Th\Ac$. Moreover $\Th\Ac$ is also quasi-free, hence its homological dimension is $\leq 1$. The results of \cite{CQ95} imply that the quotient complex $X(\Th\Ec)/F^n_{T(\Bc:\Ec)}X(\Th\Ec)$ computes the periodic cyclic homology of $\Ac$  provided that $n\geq 1$. Consequently all the complexes $F^n_{T(\Bc:\Ec)}X(\Th\Ec)$ are homotopy equivalent for $n\geq 1$. In particular taking $n=1$ one computes easily that the complex $X(\Th\Ec)/F^1_{T(\Bc:\Ec)}X(\Th\Ec)$ is isomorphic to $X(\Th\Ac)$, whence a short exact sequence of $\zz_2$-graded complexes
$$
0 \to F^1_{T(\Bc:\Ec)}X(\Th\Ec) \to X(\Th\Ec) \to X(\Th\Ac)\to 0
$$
The associated six-term cohomology exact sequence relates $H^*\big(F^1_{T(\Bc:\Ec)}X(\Th\Ec)\big)$ to $HP^*(\Ec)$ and $HP^*(\Ac)$. Now excision (\cite{CQ97}) precisely says that the natural inclusion $X(\Th\Bc)\to F^1_{T(\Bc:\Ec)}X(\Th\Ec)$ is a homotopy equivalence, which yields an isomorphism $HP^*(\Bc)\cong H^*\big(F^n_{T(\Bc:\Ec)}X(\Th\Ec)\big)$ for any $n\geq 1$. \cqfd\\

The above lemma remains unchanged if the tensor algebras $T\Ec$ and $T\Ac$ are replaced by any quasi-free extensions of $\Ec$ and $\Ac$ respectively in Diagram (\ref{univ}). Thus the cohomology groups $H^*\big(F^n_{\Rc}X(\Mch)\big)$ provide an alternative way to represent the periodic cyclic cohomology of $\Bc$. We will show in section \ref{sind} how to recover the pairing between $HP^0(\Bc)$ and the $K$-theory group $K_0(\Bc)$ in this context. For the moment observe that the morphism $E^n:H^i\big(F^n_{\Rc}X(\Mch)\big)\to HP^{i+1}(\Ac)$ factors through the universal group $H^i\big(F^n_{T(\Bc:\Ec)}X(\Th\Ec)\big)\cong HP^i(\Bc)$. We summarize these results in a corollary.
\begin{corollary}
Given any extension  $(E): 0\to\Bc\to\Ec\to\Ac\to 0$ the renormalization procedure of Proposition \ref{pren} yields a transformation 
\be
E^*\ :\ HP^i(\Bc)\to HP^{i+1}(\Ac)\ , \quad i\in \zz_2
\ee
which coincides with the connecting morphism of the extension $(E)$ given by excision.
\end{corollary}
{\it Proof:} We have only to show that the map $HP^i(\Bc)\to HP^{i+1}(\Ac)$ is the connecting morphism of the extension. Set $\Tt\Ec=\varprojlim_n T\Ec/(T(\Bc:\Ec)+J\Ec)^n$ and denote by $\iota:\Th\Ec\to\Tt\Ec$ the natural homomorphism. Also let $\pi_*:\Th\Ec\to\Th\Ac$ and $\pit_*:\Tt\Ec\to\Th\Ac$ be the homomorphisms induced by the projection $\pi:\Ec\to\Ac$. Then one has $\pit_*\circ \iota=\pi_*$, whence a commutative diagram of $\zz_2$-graded complexes and chain maps
$$
\vcenter{\xymatrix{
0 \ar[r] & F^1_{T(\Bc:\Ec)}X(\Th\Ec) \ar[r]  & X(\Th\Ec) \ar[r]^{\pi_*} \ar[d]_{\iota} & X(\Th\Ac) \ar[r]  & 0 \\
 &  & X(\Tt\Ec) \ar[ru]_{\pit_*} &   &  }}
$$
where the row is an exact sequence. Consider a linear splitting $\si:\Ac\to T\Ec$ in Diagram (\ref{univ}) as follows: first choose a linear splitting $\Ac\to\Ec$, and then map $\Ec$ into the subspace of one-tensors in $T\Ec$. The induced homomorphism $\si_*:\Th\Ac\to \Tt\Ec$ provides a right inverse for $\pit_*$: it is indeed sufficient to check the identity $\pit_*\circ\si_*=\Id_{\Th\Ac}$ on the subspace $\Ac$ which generates the whole tensor algebra $T\Ac$. \\
Then by excision, we know that any class $[\tau]\in HP^*(\Bc)$ can be represented by a cocycle $\tau$ over $F^1_{T(\Bc:\Ec)}X(\Th\Ec)$. The connecting morphism of the extension $(E)$ is nothing else but the boundary map associated to the above exact sequence of complexes: first extend $\tau$ to a linear map $\tau_R$ over $X(\Th\Ec)$. Then its coboundary $\tau_R\d$ descends to a unique cocycle $\varphi$ over $X(\Th\Ac)$ such that $\varphi\circ \pi_*=\tau_R\d$. By definition the cyclic cohomology class of $\varphi$ is the image of $[\tau]$.\\
But observe that $\varphi\circ\pit_*$ is a cocycle over $X(\Tt\Ec)$, whose pullback via the map $\iota$ is precisely $\tau_R\d$. Hence $\varphi\circ\pit_*$ is the (unique) descent of $\tau_R\d$ over $X(\Tt\Ec)$. The composite map $\tau_R\d\circ\si_* = \varphi\circ\pit_*\circ\si_*=\varphi$, which represents $E^*([\tau])$, therefore coincides with the image of $[\tau]$ under the connecting morphism. \cqfd\\

\begin{remark}\textup{Excision has been used only to show that any class in $HP^*(\Bc)$ can be represented as a cocycle over $F^n_{\Rc}X(\Mch)$ for an adequate diagram (\ref{diag}) and some $n$. Once this is known, the connecting morphism $E^*$ is given by the straightforward computation $\tau\mapsto \tau_R\d\circ\si_*$.}
\end{remark}

\section{Quasihomomorphisms}\label{squa}

The previous description of the connecting morphism is intimately related to our construction of a bivariant Chern character for quasihomomorphisms \cite{P5, P6}, if we restrict to a particular class of extensions:

\begin{definition}\label{dinv}
An extension $0\to\Bc\to\Ec\to\Ac\to 0$ is \emph{invertible} if there exists an algebra homomorphism $\rho:\Ac\to M_2(\Ec)$, such that the off-diagonal entries of the matrix $\rho$ are linear maps from $\Ac$ to the ideal $\Bc\subset\Ec$, and the upper left corner of the matrix is a linear splitting of the projection homomorphism $\Ec\to\Ac$. 
\end{definition}
From an invertible extension we construct a quasihomomorphism of odd degree as follows \cite{P5}. Let $C_1=\cc\oplus \eps\cc$ be the first Clifford algebra: it is the $\zz_2$-graded algebra generated by the unit 1 in degree zero and the element $\eps$ in degree one, with $\eps^2=1$. Define the algebra $\Ec^s=C_1\otimes\Ec^s_+ $, where $\Ec^s_+$ is the (trivially graded) matrix algebra  
\be
\Ec^s_+ = \left( \begin{matrix} \Ec & \Bc \\ \Bc & \Ec \end{matrix} \right) \subset M_2(\Ec)_ .
\ee
$\Ec^s$ is therefore $\zz_2$-graded and $\Ec^s_+$ can be identified with its subalgebra of even degree. The invertibility of the extension $0\to\Bc\to\Ec\to\Ac\to 0$ is thus equivalent to the existence of a homomorphism $\rho:\Ac\to \Ec^s_+$. Finally let $\Ic^s$ be the $\zz_2$-graded algebra $ C_1\otimes M_2(\cc)$. Then $\Ic^s\otimes\Bc=C_1\otimes M_2(\Bc)$ is a two-sided ideal in $\Ec^s$. This situation is depicted through a \emph{quasihomomorphism} of odd degree from $\Ac$ to $\Bc$:
 \be
\rho: \Ac \to \Ec^s\triangleright \Ic^s\otimes\Bc\ .
\ee
The Chern character of this quasihomomorphism lives in the bivariant cyclic cohomology of $\Ac$ and $\Bc$. The construction of \cite{P5} uses the formalism of section \ref{scm}. Thus consider a lifting $0\to\Rc\to\Mc\to\Pc\to 0$ of $0\to\Bc\to\Ec\to\Ac\to 0$, that is (Definition \ref{dlift}) a diagram of extensions 
\be
\vcenter{\xymatrix{
 & 0 \ar[d] & 0 \ar[d] & 0 \ar[d] &  \\
0 \ar [r] & \Jc \ar[r] \ar[d] & \Nc \ar[r] \ar[d] & \Qc  \ar[r] \ar[d] & 0  \\
0 \ar [r] & \Rc \ar[r] \ar[d] & \Mc \ar[r] \ar[d] & \Pc  \ar[r] \ar[d] & 0  \\
0 \ar [r] & \Bc \ar[r] \ar[d] & \Ec \ar[r] \ar[d] & \Ac  \ar[r] \ar[d] & 0  \\
 & 0 & 0 & 0 &  }} 
\ee
Notice that we do not require the extension $0\to\Rc\to\Mc\to\Pc\to 0$ be invertible. $\Jc$ and $\Rc$ are ideals respectively in $\Nc$ and $\Mc$. Moreover $\Jc=\Rc\cap\Nc$. We introduce as above the $\zz_2$-graded algebras $\Nc^s= C_1\otimes\Nc^s_+$ and $\Mc^s= C_1\otimes\Mc^s_+$, with
\be
\Nc^s_+ = \left( \begin{matrix} \Nc & \Jc \\ \Jc & \Nc \end{matrix} \right)\subset M_2(\Nc) \ ,\qquad
\Mc^s_+ = \left( \begin{matrix} \Mc & \Rc \\ \Rc & \Mc \end{matrix} \right) \subset M_2(\Mc) \ .
\ee
By construction $\Nc^s_+$ is a two-sided ideal in $\Mc^s_+$ and coincides with the kernel of the projection homomorphism $\Mc^s_+\to \Ec^s_+$. Choose a linear lifting $\ell:\Ac\to\Mc^s_+$ of the homomorphism $\rho:\Ac\to\Ec^s_+$. By the universal property of the tensor extension, one gets an homomorphism $\rho_*:T\Ac\to \Mc^s_+$ compatible with the ideals:
\be
\vcenter{\xymatrix{
0 \ar[r] & J\Ac \ar[r] \ar[d]_{\rho_*} & T\Ac \ar[r] \ar[d]_{\rho_*} & \Ac \ar[r] \ar[d]^{\rho} \ar@{.>}[dl]^{\ell} & 0 \\
0\ar[r] & \Nc^s_+ \ar[r]  & \Mc^s_+ \ar[r] & \Ec^s_+ \ar[r]  & 0 }}
\ee
Define $F=\eps\otimes \bigl( \begin{smallmatrix} 1 & 0 \\ 0 & -1 \end{smallmatrix} \bigr) $ acting on $\Mc^s$ as a multiplier of odd degree. Hence $F^2=1$ and the commutator $[F,\Mc^s_+]$ coincides with the subspace $\eps\otimes\bigl( \begin{smallmatrix} 0 & \Rc \\ \Rc & 0 \end{smallmatrix} \bigr)\subset\Mc^s$. Denote by $\tr_s$ the supertrace of odd degree $C_1\otimes M_2(\cc)\to\cc$, sending the matrix $\eps\otimes \bigl( \begin{smallmatrix} a & b \\ c & d \end{smallmatrix} \bigr)$ to $-\sqrt{2i}(a+d)$ and $1\otimes \bigl( \begin{smallmatrix} a & b \\ c & d \end{smallmatrix} \bigr)$ to 0 (see \cite{P5}). Then for any \emph{odd} integer $n$ one constructs a chain map $\chih^n$ from the $(b+B)$-complex of non-commutative differential forms over $\Mc^s_+$, to the $X$-complex of $\Mc$ as follows. $\chih^n$ has two components $\chih^n_0:\Om^n \Mc^s_+\to \Rc^n$ and $\chih^n_1:\Om^{n+1}\Mc^s_+\to \nat(\Rc^n\dd\Mc)$ given by
\be
\chih^n_0(x_0\dd x_1\ldots\dd x_n) = -\frac{\Gamma(1+\frac{n}{2})}{(n+1)!} \sum_{\la\in \La_{n+1}} \pm\, \tr_s(x_{\la(0)}[F,x_{\la(1)}]\ldots [F,x_{\la(n)}])\label{chi}
\ee
$$
\chih^n_1(x_0\dd x_1\ldots\dd x_{n+1}) = -\frac{\Gamma(1+\frac{n}{2})}{(n+1)!} \sum_{i=1}^{n+1}  \tr_s\nat(x_0[F,x_1]\ldots\dd x_i \ldots [F,x_{n+1}])
$$
where $\La_{n+1}$ is the cyclic permutation group of $n+1$ elements and $\pm$ denotes the signature of permutation $\la$. The overall minus sign is conventional (it is cancelled by the other minus coming from the supertrace). $\chih^n$ is actually defined on the direct product space $\Omh\Mc^s_+=\prod_{k\geq 0}\Om^k\Mc^s_+$ because it vanishes on differential forms of degree $> n+1$, and its image lies in the subcomplex $F^{2n-1}_{\Rc}X(\Mc)$. It clearly extends to a chain map $\chih^n:\Omh \Mch^s_+ \to F^{2n-1}_{\Rc}X(\Mch)$ where $\Mch^s_+$ is the $\Nc^s_+$-adic completion of $\Mc^s_+$. The \emph{bivariant Chern character} of degree $n$ (odd) associated to the quasihomomorphism $\rho$ is the composition of chain maps
\be
\ch^n(\rho): X(\Th\Ac) \stackrel{\gamma}{\longrightarrow} \Omh\Th\Ac \stackrel{\rho_*}{\longrightarrow} \Omh\Mch^s_+ \stackrel{\chih^n}{\longrightarrow} F^{2n-1}_{\Rc}X(\Mch) \ ,
\ee
where $\gamma$ is the generalized Goodwillie equivalence of Cuntz-Quillen (see \cite{P5}) and the middle arrow is the map of $(b+B)$-complexes induced by the homomorphism $\rho_*:\Th\Ac\to\Mch^s_+$. Hence if $\tau$ is a cocycle over $F^{2n-1}_{\Rc}X(\Mch)$, the composite $\tau \ch^n(\rho)$ defines a periodic cyclic cohomology class over $\Ac$.

\begin{proposition}\label{pbiv}
Let $(E):0\to\Bc\to\Ec\to\Ac\to 0$ be an invertible extension with lifting $0\to\Rc\to\Mc\to\Pc\to 0$, and let $\rho:\Ac\to \Ec^s\triangleright\Ic^s\otimes\Bc$ be the associated quasihomomorphism. Then for any cocycle $\tau$ over $F^{2n-1}_{\Rc}X(\Mch)$ representing a class $[\tau]\in HP^i(\Bc)$, where $n$ is odd, the equality 
\be
\tau \ch^n(\rho) = \sqrt{2\pi i} \, E^*([\tau])
\ee
holds in $HP^{i+1}(\Ac)$.
\end{proposition}
{\it Proof:} We shall relate $\tau\ch^n(\rho)$ to $E^*([\tau])$ via the eta-cochain and its renormalization introduced in \cite{P5, P6}. The eta-cochain of degree $n+1$ has two components $\etah^{n+1}_0:\Om^{n+1}\Mc^s_+\to \Rc^{n+1}$ and $\etah^{n+1}_1:\Om^{n+2}\Mc^s_+\to \nat(\Rc^{n+1}\dd\Mc)$ given by
\beq
\lefteqn{\etah^{n+1}_0(x_0\dd x_1\ldots\dd x_{n+1}) =  \frac{\Gamma(\frac{n}{2}+1)}{(n+2)!} \, \frac{1}{2}\tr_s\Big (F x_0[F,x_1]\ldots [F,x_{n+1}]+ } \non \\
&&\qquad \qquad \qquad  \sum_{i=1}^{n+1}(-)^{(n+1)i} [F,x_i]\ldots [F,x_{n+1}] Fx_0 [F,x_1]\ldots [F,x_{i-1}] \Big) \non
\eeq
\beq
\lefteqn{\etah^{n+1}_1(x_0\dd x_1\ldots\dd x_{n+2}) =} \non\\
&& \frac{\Gamma(\frac{n}{2}+1)}{(n+3)!} \sum_{i=1}^{n+2}  \frac{1}{2}\tr_s\nat(ix_0 F + (n+3-i)Fx_0)[F,x_1]\ldots\dd x_i \ldots [F,x_{n+2}]\ .\non
\eeq
These components extend as above to a linear map $\etah^{n+1}:\Omh\Mch^s_+\to F^{2n+1}_{\Rc}X(\Mch)$. The eta-cochain makes the connection between the chain maps $\chih^n$ and $\chih^{n+2}$ for any odd integer $n$. Indeed let $\d$ and $(b+B)$ denote the boundaries on the complexes $X(\Mch)$ and $\Omh\Mch^s_+$ respectively. The following transgression relation holds (\cite{P5}):
\be
\chih^n - \chih^{n+2}= \d\circ \etah^{n+1} - \etah^{n+1}\circ(b+B)\ .\label{trans}
\ee
Now let $\tau$ be a cocycle over the subcomplex $F^{2n-1}_{\Rc}X(\Mch)$, and choose any renormalization $\tau_R:X(\Mch)\to \cc$ as in section \ref{scm}. Hence the composite map $\tau_R\d$ vanishes on $F^{2n-1}_{\Rc}X(\Mch)$ but not on $X(\Mch)$. Define the cochain $\chi_R:\Omh\Mch^s_+\to\cc$ by
$$
\chi_R :=  \sum_{k\ \mathrm{odd}\ <n}-\tau_R\d\etah^{k+1}\ .
$$
Then $\chi_R$ is a $(b+B)$-cocycle cohomologous to $\tau\chih^n$. Indeed using the transgressions (\ref{trans}) one gets
$$
\tau\chih^n - \chi_R = \Big(\sum_{k\ \mathrm{odd}\ <n} \tau_R\etah^{k+1}\Big)\circ(b+B)\ .
$$
Hence composition with the chain map $\rho_*\gamma: X(\Th\Ac)\to \Omh\Mch^s_+$ yields the equality of cyclic cohomology classes $\tau\ch^n(\rho)\equiv \chi_R \rho_*\gamma$ in $HP^{i+1}(\Ac)$. It remains to compare $\chi_R \rho_*\gamma$ and $E^*([\tau])$. We produce a deformation of the homomorphism $\rho_*:T\Ac\to\Mc^s_+$ as follows. Recall that $\rho_*$ is induced by a linear lifting $\ell:\Ac\to\Mc^s_+$ of the homomorphism $\rho:\Ac\to\Ec^s_+$. In matrix form we can write
$$
\ell = \left( \begin{matrix} \si & \la \\ \mu & \sib \end{matrix} \right)\ ,
$$
where $\si:\Ac\to\Mc$ is a linear splitting of the projection homomorphism $\Mc\to\Ac$ as in section \ref{scm}, and $\la,\mu$ are linear maps from $\Ac$ to $\Rc$. Consider the linear homotopy of linear maps $\ell^t:\Ac\to \Mc^s_+$, defined for any $t\in [0,1]$ by
$$
\ell^t = \left( \begin{matrix} \si & t\la \\ t\mu & \sib \end{matrix} \right)\ .
$$
In particular $\ell^0= \bigl( \begin{smallmatrix} \si & 0 \\ 0 & \sib \end{smallmatrix} \bigr)$ is a diagonal matrix, and $\ell^1=\ell$. Then observe that $\ell^t$ followed by the projection $\Mc^s_+= \bigl( \begin{smallmatrix} \Mc & \Rc \\ \Rc & \Mc \end{smallmatrix} \bigr) \to  \bigl( \begin{smallmatrix} \Ac & 0 \\ 0 & \Ac \end{smallmatrix} \bigr)=\Ac\oplus\Ac$ yields an algebra homomorphism (independent of $t$). Hence one gets as usual an homomorphism $\rho_*^t:T\Ac\to \Mc^s_+$ by means of the diagram
$$
\vcenter{\xymatrix{
0 \ar[r] & J\Ac \ar[r] \ar[d]_{\rho_*^t} & T\Ac \ar[r] \ar[d]_{\rho_*^t} & \Ac \ar[r] \ar[d] \ar@{.>}[dl]^{\ell^t} & 0 \\
0\ar[r] & M_2(\Rc) + \Nc^s_+   \ar[r]  & \Mc^s_+ \ar[r] & \Ac \oplus \Ac \ar[r]  & 0 }}
$$
where the ideal $ M_2(\Rc) + \Nc^s_+ = \bigl( \begin{smallmatrix} \Rc+\Nc & \Rc \\ \Rc & \Rc+\Nc \end{smallmatrix} \bigr)$ is the kernel of the projection $\Mc^s_+\to \Ac\oplus\Ac$. Thus $\rho_*^t$ extends to an homomorphism from $\Th\Ac$ to the pro-algebra $\Mct^s_+ = \varprojlim_m \Mc^s_+/(M_2(\Rc)+\Nc^s_+)^m$ for any $t\in [0,1]$. This provides a homotopy between $\rho_*^0= \bigl( \begin{smallmatrix} \si_* & 0 \\ 0 & \sib_* \end{smallmatrix} \bigr)$ and $\rho_*^1=\rho_*$. One knows that $\tau_R\d$ vanishes on $F^{2n-1}_{\Rc}X(\Mch)$. Hence $\chi_R$ descends to a $(b+B)$-cocycle over $\Omh\Mct^s_+$, and by homotopy invariance of periodic cyclic cohomology \cite{CQ95}, the composite maps $\chi_R\rho_*^t\gamma: X(\Th\Ac)\to\cc$ define the same cohomology class for all $t$. We deduce the equality of cyclic cohomology classes
$$
\tau\ch^n(\rho) \equiv \chi_R\rho_*\gamma=\chi_R\rho_*^1\gamma \equiv \chi_R\rho_*^0\gamma=\sum_{k\ \mathrm{odd}\ <n}-\tau_R\d\etah^{k+1}\rho_*^0\gamma 
$$
in $HP^{i+1}(\Ac)$. Since $\rho_*^0$ is a diagonal matrix, $[F,\rho_*^0(x)]=0$ for any $x\in \Th\Ac$ and the map $\etah^{k+1}\rho_*^0$ vanishes unless $k=-1$. Hence only the term containing $\etah^0\rho_*^0$ survives in the sum over $k$. One has
\beq
\eta^0_0\rho_*^0(x) &=& \Gamma(1/2)\, \frac{1}{2}\tr_s(F\rho_*^0(x))\ =\ -\sqrt{2\pi i}\, \frac{1}{2}(\si_*(x)-\sib_*(x))\ , \non\\
\eta^0_1\rho_*^0(x\dd y) &=& \Gamma(1/2)\, \frac{1}{2}\tr_s\nat(F\rho_*^0(x)\dd\rho_*^0(y)) \non\\
&=& -\sqrt{2\pi i}\, \frac{1}{2}\nat(\si_*(x)\dd\si_*(y)-\sib_*(x)\dd\sib_*(y))\ , \non
\eeq
which shows that $\tau\ch^n(\rho)$ is cohomologous to $\sqrt{2\pi i}\,\tau_R\d\circ \frac{1}{2}(\si_*-\sib_*)$, where $\si_*$ and $\sib_*$ are viewed as chain maps from $X(\Th\Ac)$ to $X(\Mct)$. Finally consider $\rho_*^t$ as a family of homomorphisms $\Th\Ac\to M_2(\Mct)$. The cup-product of the cocycle $\tau_R\d$ over $X(\Mct)$ with the usual trace over $M_2(\cc)$ yields a cocycle $\tau_R\d\#\tr$ over $X(M_2(\Mct))$. Using homotopy invariance, one has the equality of periodic cyclic cohomology classes 
$$
\tau_R\d\si_* + \tau_R\d\sib_* = (\tau_R\d\#\tr)\rho_*^0 \equiv (\tau_R\d\#\tr)\rho_*^1 = (\tau_R\d\#\tr)\rho_* \equiv 0
$$
where the last equality comes from Proposition \ref{pren} and the fact that $\rho:\Ac\to M_2(\Ec)$ is an homomorphism. Hence $\tau_R\d\circ \frac{1}{2}(\si_*-\sib_*)\equiv \tau_R\d\si_*=E^*([\tau])$, and $\tau\ch^n(\rho)$ coincides with $\sqrt{2\pi i}\, E^*([\tau])$ in $HP^{i+1}(\Ac)$ as wanted. \cqfd\\

\begin{remark}\textup{The above proof relates the connecting morphism of an invertible extension to the boundary of the renormalized eta-cochain, a technique introduced in \cite{P6} as a way of building \emph{local} representatives for the bivariant Chern character of quasihomomorphisms. It was shown that this procedure is intimately related to anomalies in quantum field theory. Thus we may consider the cyclic cocycle $\tau_R\d\si_*$ representing $E^*([\tau])$ as a kind of ``anomaly formula'' adapted to extensions.}
\end{remark}

\section{Index theorem}\label{sind}

First recall the definition of the algebraic $K$-theory groups of a (non-unital) algebra $\Ac$ in low degrees. As usual denote by $M_{\infty}(\Ac)=\varinjlim_NM_N(\Ac)$ the inductive limit of matrix algebras with entries in $\Ac$, under the inclusion maps $a\mapsto \bigl( \begin{smallmatrix} a & 0 \\ 0 & 0 \end{smallmatrix} \bigr)$. Let $\Ac^+=\Ac\oplus\cc$ be the unitalization of $\Ac$. An equivalence relation is defined on the set of idempotents in $M_{\infty}(\Ac^+)$ as follows: two idempotents $e$ and $e'$ are equivalent if there exists an invertible matrix $g$ with entries in $\Ac^+$, such that $g- 1\in M_{\infty}(\Ac^+)$ and $e'=g^{-1}eg$ (similarity). The set of equivalence classes of idempotents forms a semigroup for the direct sum of matrices. Denote by $K_{00}(\Ac)$ its Grothendieck group and let $K_0(\Ac)$ be the kernel of the morphism $K_{00}(\Ac^+)\to K_{00}(\cc)=\zz$. The elements of $K_0(\Ac)$ are represented by formal differences $[e]-[f]$ of equivalence classes of idempotents $e,f\in M_{\infty}(\Ac^+)$ such that $e\equiv f$ mod $ M_{\infty}(\Ac)$. Any such class can be further reduced to a difference $[e]-[p_N]$ where $p_N$ is the diagonal matrix $\bigl( \begin{smallmatrix} 1_N & 0 \\ 0 & 0 \end{smallmatrix} \bigr)$ with $N$ units on the diagonal. \\
For any integer $N$ denote by $GL_N(\Ac)$ the group of invertible matrices $g\in M_N(\Ac^+)$ such that $g\equiv 1_N$ mod $M_N(\Ac)$, and by $GL_{\infty}(\Ac)$ the inductive limit of the groups $GL_N(\Ac)$ under the inclusions $a\mapsto \bigl( \begin{smallmatrix} a & 0 \\ 0 & 1 \end{smallmatrix} \bigr)$. Then $K_1(\Ac)$ is the abelianization of $GL_{\infty}(\Ac)$, i.e. its quotient by the commutator subgroup $[GL_{\infty}(\Ac),GL_{\infty}(\Ac)]$.\\

Now let $(E): 0\to \Bc\to\Ec\to\Ac\to 0$ be an extension. The connecting map for the algebraic $K$-theory in low degree is constructed as follows \cite{M}. Take a class $[g]\in K_1(\Ac)$ represented by an invertible matrix $g\in M_N(\Ac^+)$ with $g-1\in M_N(\Ac)$. Then $\bigl( \begin{smallmatrix} g & 0 \\ 0 & g^{-1} \end{smallmatrix} \bigr)\in M_{2N}(\Ac^+)$ can be lifted to an invertible matrix in $M_{2N}(\Ec^+)$. Indeed, choose any preimages $Q,P\in M_N(\Ec^+)$ of $g,g^{-1}$ respectively. Then 
\be
G=\left( \begin{matrix} 0 & 1 \\ -1 & P \end{matrix} \right)\left( \begin{matrix} 1 & 0 \\ Q & 1 \end{matrix} \right)\left( \begin{matrix} 1 & -P \\ 0 & 1 \end{matrix} \right) = \left( \begin{matrix} Q & 1-QP \\ PQ-1 & P+P(1-QP) \end{matrix} \right)
\ee
is an invertible matrix in $M_{2N}(\Ec^+)$. Moreover the image of $1-QP$ vanishes in $M_N(\Ac^+)$, hence $1-QP\in M_N(\Bc)$ and similarly for $PQ-1$. Thus $G$ is an invertible lifting of $\bigl( \begin{smallmatrix} g & 0 \\ 0 & g^{-1} \end{smallmatrix} \bigr)$. Viewing $p_N=\bigl( \begin{smallmatrix} 1_N & 0 \\ 0 & 0 \end{smallmatrix} \bigr)$ as a $2N\times 2N$ matrix, the idempotent $e=G^{-1}p_NG$ fulfills the property $e\equiv p_N$ mod $M_{2N}(\Bc)$ and consequently $[e]-[p_N]$ determines a class in $K_0(\Bc)$. The latter is independent of the choice of lifting $G$. Since by definition any additive map from $GL_{\infty}(\Ac)$ to an abelian group factors through $K_1(\Ac)$, we obtain the \emph{index morphism} of the extension $(E)$
\be
\Ind_E\ :\ K_1(\Ac)\to K_0(\Bc)\ .
\ee

Before showing that the index map is adjoint to the connecting morphism in periodic cyclic cohomology $E^*: HP^0(\Bc)\to HP^1(\Ac)$, let us explain how Connes' pairing $HP^0(\Bc)\times K_0(\Bc)\to \cc$ is computed when the cyclic cohomology of $\Bc$ is represented by cocycles over the complexes $F^n_{\Rc}X(\Mch)$ as in section \ref{scm}. We will use the formulas established in \cite{P5} \S 4 in connection with the Chern-Connes character in cyclic homology. Consider a lifting $0\to\Rc\to\Mc\to\Pc\to 0$ of the extension $(E)$:
\be
\vcenter{\xymatrix{
 & 0 \ar[d] & 0 \ar[d] & 0 \ar[d] &  \\
0 \ar [r] & \Jc \ar[r] \ar[d] & \Nc \ar[r] \ar[d] & \Qc  \ar[r] \ar[d] & 0  \\
0 \ar [r] & \Rc \ar[r] \ar[d] & \Mc \ar[r] \ar[d] & \Pc  \ar[r] \ar[d] & 0  \\
0 \ar [r] & \Bc \ar[r] \ar[d] & \Ec \ar[r] \ar[d] & \Ac  \ar[r] \ar[d] & 0  \\
 & 0 & 0 & 0 &  }} \label{diag2}
\ee 
For any $n\geq 1$, the intersection $\Rc\cap\Nc^n$ is an ideal in $\Rc$ and the quotient $\Rc/(\Rc\cap\Nc^n)$ is a nilpotent extension of $\Bc$. Take $\Rch$ as the projective limit $\varprojlim_n \Rc/(\Rc\cap\Nc^n)$. Proceeding as in \cite{CQ95}, any idempotent $e\in M_{\infty}(\Bc^+)$ such that $e-p_N\in M_{\infty}(\Bc)$ can be lifted to an idempotent $\eh\in M_{\infty}(\Rch^+)$. The latter is defined up to similarity in the matrix algebra over $\Rch^+$. One has $\eh-p_N\in M_{\infty}(\Rch)$ and the trace $\tr(\eh-p_N)\in \Rch$ defines a cycle of even degree in the subcomplex $F^1_{\Rc}X(\Mch)\subset X(\Mch)$. Thus if $\tau$ is a cocycle of even degree over $F^1_{\Rc}X(\Mch)$, the pairing $\tau\#\tr(\eh-p_N)$ is defined. More generally, $\tr\big((\eh-p_N)^{2n+1}\big)$ is a cycle of even degree in $F^{4n+1}_{\Rc}X(\Mch)$ and hence can be paired with any cocycle $\tau:F^{4n+1}_{\Rc}X(\Mch)\to \cc$.
\begin{lemma}\label{lpair}
Let $\tau:F^{4n+1}_{\Rc}X(\Mch)\to\cc$ be a cocycle of even degree. Let $e\in M_{\infty}(\Bc^+)$ be an idempotent such that $e-p_N\in M_{\infty}(\Bc)$, and choose an idempotent lifting $\eh\in M_{\infty}(\Rch^+)$ of $e$. Then the formula
\be
\langle [\tau], [e] \rangle = \tau\#\tr \big( (\eh-p_N)^{2n+1}\big) \label{pair}
\ee
descends to a well-defined pairing $\varinjlim_nH^0\big(F^{4n+1}_{\Rc}X(\Mch)\big)\times K_0(\Bc)\to\cc$.
\end{lemma} 
{\it Proof:} We have to show that $\tau\#\tr \big( (\eh-p_N)^{2n+1}\big)$ does not depend on $n$ (sufficiently large), and that it is invariant when $\eh$ is conjugated by an invertible matrix $u$ such that $u-1\in M_{\infty}(\Rch^+)$ and $u^{-1}\eh u\equiv p_N$ mod $M_{\infty}(\Rch)$. For convenience we rewrite the pairing using the $\zz_2$-graded algebra of $2\times 2$ matrices over $M_{\infty}(\Rch^+)$, with grading induced by the decomposition of matrices into diagonal/off diagonal form: consider the odd element $F=\bigl( \begin{smallmatrix} 0 & 1 \\ 1 & 0 \end{smallmatrix} \bigr)$ such that $F^2=1$, and set $f=\bigl( \begin{smallmatrix} \eh & 0 \\ 0 & p_N \end{smallmatrix} \bigr)$ as an idempotent of even degree. Then if $\tr_s$ denotes the supertrace on $M_2(M_{\infty}(\cc))$ one has
$$
\tau\#\tr \big( (\eh-p_N)^{2n+1}\big) = \tau\#\tr_s \big( F([F,f])^{2n+1}\big)\ .
$$
The right-hand-side is recognized as a Chern-Connes pairing \cite{C86} and has well-known properties. In particular it does not depend on $n$ provided $F([F,f])^{2n+1}$ remains in the domain of the supertrace $\tau\#\tr_s$, and it is invariant with respect to homotopies of $f$ preserving the condition $[F, f]\in M_2(M_{\infty}(\Rch))$. Now let $u$ be an invertible matrix such that $u-1\in M_{\infty}(\Rch^+)$ and $u^{-1}\eh u\equiv p_N$ mod $M_{\infty}(\Rch)$. Let $v$ be the image of $u$ under the projection $\Rch^+=\Rch\oplus\cc\to\cc$. Then $v$ is an invertible matrix such that $v-1\in  M_{\infty}(\cc)$ and $v^{-1}p_Nv=p_N$. The invertible matrix of even degree $g=\bigl( \begin{smallmatrix} u & 0 \\ 0 & v \end{smallmatrix} \bigr)$ conjugates $f$ to $g^{-1}fg=\bigl( \begin{smallmatrix} u^{-1}\eh u & 0 \\ 0 & p_N \end{smallmatrix} \bigr)$, and fulfills the commutation relation $[F,g]\in M_2(M_{\infty}(\Rch))$. This allows to construct a (stable) homotopy between $f$ and $g^{-1}fg$ by a standard procedure using rotation matrices. \cqfd\\

\begin{remark}\textup{In general the definition of $\Rch$ given here does not coincide with the pro-algebra $\varprojlim_n\Rc/\Jc^n$. It does coincide under strong conditions, for example when the equality $\Rc\cap\Nc^n=\Jc^n$ holds for all $n$. The latter condition was implicitely assumed in \cite{P5}, where the cycle $\tr \big( (\eh-p_N)^{2n+1}\big)$ was taken as the definition of the Chern character in $K$-theory.}
\end{remark}

Finally recall Connes' pairing $HP^1(\Ac)\times K_1(\Ac)\to \cc$ in the Cuntz-Quillen formalism \cite{CQ95}. Let $[\varphi]\in HP^1(\Ac)$ be a cyclic cohomology class represented by a cocycle of odd degree $\varphi:X(\Th\Ac)\to\cc$, where $\Th\Ac$ is the $J\Ac$-adic completion of the tensor algebra $T\Ac$. Let $[g]\in K_1(\Ac)$ be represented by an invertible element $g\in GL_{\infty}(\Ac)$. Then $g$ can be lifted to an invertible element $\gh\in GL_{\infty}(\Th\Ac)$, and the one-form $\nat(\gh^{-1}\dd\gh)\in \Om^1\Th\Ac_{\nat}$ is a cycle of odd degree in the complex $X(\Th\Ac)$ whose homology class is independent of the choice of lifting. The pairing is defined as
\be
\langle [\varphi], [g] \rangle = \frac{1}{\sqrt{2\pi i}}\, \varphi \big(\gh^{-1}\dd\gh\big)\ .
\ee
One can think of the normalization factor $1/\sqrt{2\pi i}$ as a pure convention. However note that this normalization is uniquely determined by the compatibility of the bivariant Chern character with Bott periodicity, see for example \cite{P5}.

\begin{theorem}\label{t}
Let $(E):0\to\Bc\to\Ec\to\Ac\to 0$ be an extension with lifting $0\to\Rc\to\Mc\to\Pc\to 0$. Let $\tau:F^{4n+1}_{\Rc}X(\Mch)\to\cc$ ce a cocycle representing an element $[\tau]\in HP^0(\Bc)$, and take any $[g]\in K_1(\Ac)$. Then 
\be
\langle [\tau], \Ind_E([g]) \rangle = \sqrt{2\pi i}\, \langle E^*([\tau]), [g] \rangle\ .
\ee
\end{theorem}
{\it Proof:} We will not use excision since we assume from the beginning that the cyclic cohomology class $[\tau]$ is represented by a cocycle over $F^{4n+1}_{\Rc}X(\Mch)$. Let $g\in GL_N(\Ac)$ be an invertible element representing $[g]$. Thus in particular $g-1_N\in M_N(\Ac)$. We shall replace $(E)$ with an \emph{invertible} extension as follows. Denote by $\cc[z,z^{-1}]$ the commutative algebra of Laurent polynomials in the indeterminate variable $z$, and let $\Cc$ be the subalgebra of polynomials $f\in \cc[z,z^{-1}]$ such that $f(1)=0$. Equivalently, $\Cc$ is the (non-unital) commutative algebra generated by two elements $u,v$ with relations $uv=vu=-u-v$. The inclusion of $\Cc$ into $\cc[z,z^{-1}]$ is recovered by setting $z=1+u$ and $z^{-1}=1+v$. The geometric picture is that of the algebra of trigonometric functions over the unit circle, vanishing at point $z=1$. Hence $\Cc$ is a suitable algebraic definition of a suspension algebra. We define a homomorphism
$$
\al: \Cc\to M_N(\Ac)
$$
by setting $\al(u)=g-1$ and $\al(v)=g^{-1}-1$. Equivalently we may extend $\al$ to a unital homomorphism from $\Cc^+=\cc[z,z^{-1}]$ to $M_N(\Ac^+)$ and set $\al(z)=g$, $\al(z^{-1})=g^{-1}$. Thus $\al$ carries the ``Bott element'' $[z]\in K_1(\Cc)$ to $[g]\in K_1(\Ac)$. Define $(F)$ as the \emph{pullback extension} of $(E)$ (tensored with $M_N(\cc)$) induced by $\al$, that is, $(F)$ is the first row in the commutative diagram 
$$
\vcenter{\xymatrix{
0 \ar[r] & M_N(\Bc) \ar[r] \ar@{=}[d] & \Fc \ar[r] \ar[d]_{\beta} & \Cc \ar[r] \ar[d]^{\al} & 0 \\
0\ar[r] & M_N(\Bc) \ar[r]  & M_N(\Ec) \ar[r]^{\pi} & M_N(\Ac) \ar[r]  & 0 }}
$$
Explicitely $\Fc=\{(h,f)\in M_N(\Ec)\oplus \Cc \,|\, \pi(h)=\al(f)\}$. The homomorphisms $\Fc\to M_N(\Ec)$ and $\Fc\to\Cc$ are induced respectively by the projections onto the first and second summand in $M_N(\Ec)\oplus \Cc$. It is immediate from the construction of the index map that the equality
$$
\Ind_E([g]) = \Ind_F([z])
$$
holds in $K_0(\Bc)$. Now use the universal property of the tensor algebra to extend $\al$ to a morphism from the extension $0\to J\Cc\to T\Cc\to \Cc\to 0$ to the third column of Diagram (\ref{diag2}) tensored with $M_N(\cc)$. This yields a homomorphism $\al_*: T\Cc \to M_N(\Pc)$ respecting the ideals $J\Cc\subset T\Cc$ and $M_N(\Qc)\subset M_N(\Pc)$. Again we define the pullback extension
$$
\vcenter{\xymatrix{
0 \ar[r] & M_N(\Rc) \ar[r] \ar@{=}[d] & \Gc \ar[r] \ar[d]_{\beta_*} & T\Cc \ar[r] \ar[d]^{\al_*} & 0 \\
0\ar[r] & M_N(\Rc) \ar[r]  & M_N(\Mc) \ar[r] & M_N(\Pc) \ar[r]  & 0 }}
$$
and similarly the homomorphism $\al_*:J\Cc\to M_N(\Qc)$ restricted to the ideals yields a pullback $0\to M_N(\Jc) \to \Hc\to T\Cc\to 0$. All these extensions fit together in a commutative diagram
$$
\vcenter{\xymatrix{
 & 0 \ar[d] & 0 \ar[d] & 0 \ar[d] &  \\
0 \ar [r] & M_N(\Jc) \ar[r] \ar[d] & \Hc \ar[r] \ar[d] & J\Cc  \ar[r] \ar[d] & 0  \\
0 \ar [r] & M_N(\Rc) \ar[r] \ar[d] & \Gc \ar[r] \ar[d] & T\Cc  \ar[r] \ar[d] & 0  \\
0 \ar [r] & M_N(\Bc) \ar[r] \ar[d] & \Fc \ar[r] \ar[d] & \Cc  \ar[r] \ar[d] & 0  \\
 & 0 & 0 & 0 &  }}
$$
This diagram is naturally mapped to (\ref{diag2}) (tensored with $M_N(\cc)$) under the homomorphism $\beta_*:\Gc\to M_N(\Mc)$. If $\tau$ is any cocycle over $F^{4n+1}_{\Rc}X(\Mch)$, its cup-product with the trace of matrices yields a cocycle $\tau\#\tr$ over $F^{4n+1}_{M_N(\Rc)}X(M_N(\Mch))$, which may be pulled back to a cocycle $\beta^*(\tau)$ over $F^{4n+1}_{M_N(\Rc)}X(\Gch)$. One has
$$
\langle [\tau], \Ind_E([g]) \rangle = \langle \beta^*([\tau]), \Ind_F([z]) \rangle \ .
$$
Moreover the homomorphism $\al: \Cc\to\Ac$ induces a pullback in cyclic cohomology $\al^*: HP^*(\Ac)\to HP^*(\Cc)$. The pair $\al^*,\beta^*$ intertwines the action of the connecting morphisms $E^*$ and $F^*$ in the sense that $\al^*\circ E^* = F^*\circ \beta^*$. Therefore one has 
$$
\langle E^*([\tau]), [g] \rangle = \langle E^*([\tau]), \al([z]) \rangle = \langle \al^*\circ E^*([\tau]), [z] \rangle = \langle F^*\circ \beta^*([\tau]), [z] \rangle
$$
and the equality $\langle [\tau], \Ind_E([g]) \rangle = \sqrt{2\pi i}\, \langle E^*([\tau]), [g] \rangle$ would follow from the equality $\langle [\tau'], \Ind_F([z]) \rangle = \sqrt{2\pi i}\, \langle F^*([\tau']), [z] \rangle$ for the cocycle $\tau'=\beta^*(\tau)$. We decided to replace the extension $(E)$ with the extension $(F)$  because the latter is invertible. Indeed choose arbitrary liftings $U,V\in \Fc$ of $u,v\in\Cc$ and set $Q=1+U$, $P=1+V$ in $\Fc^+$. Then $1-QP$ and $1-PQ$ sit in the ideal $M_N(\Bc)$, and the map $\rho$ defined on generators by
$$
\rho(u) = \left( \begin{matrix} U & 1-QP \\ PQ-1 & V+P(1-QP) \end{matrix} \right)\ ,\quad \rho(v) = \left( \begin{matrix} V+P(1-QP) & PQ-1 \\ 1-QP & U \end{matrix} \right) 
$$
extends to a homomorphism $\rho:\Cc\to \Fc^s_+$. Passing to the unitalized algebra $\Cc^+$, the map $\rho$ carries $z$ to the invertible $\bigl( \begin{smallmatrix} Q & 1-QP \\ PQ-1 & P+P(1-QP) \end{smallmatrix} \bigr)$ and $z^{-1}$ to its inverse $\bigl( \begin{smallmatrix} P+P(1-QP) & PQ-1 \\ 1-QP & Q \end{smallmatrix} \bigr)$. Observe that the $K$-theory class $\Ind_F([z])$ is represented by the idempotent $e=\rho(z)^{-1}\bigl( \begin{smallmatrix} 1 & 0 \\ 0 & 0 \end{smallmatrix} \bigr)\rho(z) \in M_2(M_N(\Bc)^+)$. As explained in section \ref{squa} the invertible extension $(F)$ determines a quasihomomorphism 
$$
\rho:\Cc\to \Fc^s \triangleright \Ic^s\otimes M_N(\Bc)\ .
$$
Its bivariant Chern character  $\ch^{2n+1}(\rho)$ is a chain map $X(\Th\Cc) \to F^{4n+1}_{M_N(\Rc)}X(\Gch)$. By Proposition \ref{pbiv}, the composite $\tau'\ch^{2n+1}(\rho)$ represents a cyclic cohomology class over $\Cc$ which coincides with $\sqrt{2\pi i}\, F^*([\tau'])$. Let us calculate explicitely the pairing $\langle \tau'\ch^{2n+1}(\rho), [z] \rangle$. We know that $\rho$ lifts to a homomorphism $\rho_*:\Th\Cc\to \Gch^s_+=\varprojlim_k \Gc^s_+/(\Hc^s_+)^k$. Choose an invertible lifting $\zh\in (\Th\Cc)^+$ of $z$. The idempotent 
$$
\eh = \rho_*(\zh)^{-1} \left( \begin{matrix} 1 & 0 \\ 0 & 0 \end{matrix} \right) \rho_*(\zh) \ \in M_2\Big(\widehat{M_N(\Rc)}^+\Big)\ ,
$$
where $\widehat{M_N(\Rc)} = \varprojlim_k M_N(\Rc)/(M_N(\Rc)\cap \Hc^k)$, is a lifting of the idempotent $e=\Ind_F(z)$. According to Lemma \ref{lpair} the pairing $\langle [\tau'], \Ind_F([z]) \rangle$ is given by $\tau'\#\tr \big( (\eh-p_1)^{2n+1}\big)$. On the other hand, the calculation performed in the proof of \cite{P5} Theorem 6.3 part III) applies verbatim and yields
$$
\langle \tau'\ch^{2n+1}(\rho), [z] \rangle = \tau'\#\tr \big( (\eh-p_1)^{2n+1}\big)=\langle [\tau'], \Ind_F([z]) \rangle
$$
where $p_1$ is the $2\times 2$ matrix $\bigl( \begin{smallmatrix} 1 & 0 \\ 0 & 0 \end{smallmatrix} \bigr)$. Hence we conclude that $\langle [\tau'], \Ind_F([z]) \rangle = \sqrt{2\pi i}\,\langle F^*([\tau'], [z] \rangle$ as wanted. \hfill\cqfd\\

\begin{corollary}
Let $(E):0\to\Bc\to\Ec\to\Ac\to 0$ be an extension with lifting $0\to\Rc\to\Mc\to\Pc\to 0$. Then for any class $[\tau]\in HP^0(\Bc)$ represented by a cocycle of even degree $\tau: F^{4n+1}_{\Rc}X(\Mch)\to\cc$ and any class $[g]\in K_1(\Ac)$ represented by an invertible $g\in GL_{\infty}(\Ac)$, one has
\be
\langle [\tau], \Ind_E([g]) \rangle = \tau_R\d\#\tr \big( \gt^{-1}\dd\gt \big)\ .
\ee
Here $\tau_R: X(\Mch)\to\cc$ is any renormalization of $\tau$, and $\gt\in GL_{\infty}(\Mct)$ is any invertible lifting of $g$ in the pro-algebra $\Mct=\varprojlim_k\Mc/(\Rc+\Nc)^k$.
\end{corollary}
{\it Proof:} Choose a linear section $\si:\Ac\to\Mc$ of the projection homomorphism. By definition the cyclic cohomology class $E^*([\tau]) \in HP^1(\Ac)$ is represented by $\tau_R\d\circ \si_* : X(\Th\Ac) \to X(\Mct) \to \cc$, where $\si_*$ is the homomorphism $\Th\Ac\to \Mct$ induced by $\si$ (section \ref{scm}). Hence by Theorem \ref{t}, evaluating this cocycle on $[g]$ yields the formula
$$
\langle [\tau], \Ind_E([g]) \rangle = \tau_R\d\#\tr \big( \si_*(\gh)^{-1}\dd\si_*(\gh) \big)
$$
for any invertible lifting $\gh\in GL_{\infty}(\Th\Ac)$. Then by \cite{CQ95} \S 12 the homology class of $\nat \si_*(\gh)^{-1}\dd\si_*(\gh)$ remains unchanged if $\si_*(\gh)$ is replaced by any other invertible lifting $\gt\in GL_{\infty}(\Mct)$ of $g$.  \cqfd\\

\section{Pseudodifferential operators}\label{spseu}

Let $\pi: M\to B$ be a proper submersion of smooth manifolds without boundary. Hence at any point $b\in B$ the fiber $\pi^{-1}(b)=(M/B)_b$ is a compact submanifold of $M$. Define $E\to B$ as the infinite-dimensional Fr\'echet bundle whose fiber at a point $b$ is the space of smooth functions $\cinf((M/B)_b)$. The space of smooth sections $\cinf(B,E)$ is thus isomorphic to $\cinf(M)$. For any $k\in\zz$ we denote by $CL^k_b=CL^k((M/B)_b)$ the space of \emph{classical} pseudodifferential operators of order $k$ on the manifold $(M/B)_b$. In particular $CL^{-1}_b$ is a two-sided ideal in the algebra $CL^0_b$ and the quotient $CL^0_b/CL^{-1}_b=LS^0_b$ is isomorphic to the commutative algebra of smooth functions $\cinf(S^*(M/B)_b)$ on the cotangent sphere bundle of $(M/B)_b$. The projection homomorphism $CL^0_b\to LS^0_b$ is the map which carries a pseudodifferential operator of order zero to its leading symbol. The algebra of smooth families of fiberwise pseudodifferential operators $\cinf(B,CL^0)$, parametrized by the base manifold $B$, naturally acts on $\cinf(E)$ by endomorphisms. The following algebras of smooth families of operators \emph{with compact support} on $B$
\be
\Bc= \cinfc(B,CL^{-1}) \ ,\quad \Ec=\cinfc(B,CL^0) \ ,\quad \Ac=\cinfc(B,LS^0) 
\ee 
thus lead to an extension $(E):0\to\Bc\to\Ec\to\Ac\to 0$. The index morphism $\Ind_E: K_1(\Ac)\to K_0(\Bc)$ maps a family of elliptic symbols $g\in GL_{\infty}(\Ac)$ to an idempotent in $M_{\infty}(\Bc^+)$ representing a $K$-theory class (index bundle) of the base manifold $B$. Our aim is to evaluate the image of $\Ind_E$ on certain cyclic cohomology classes $[\tau]\in HP^0(\Bc)$ associated to closed currents over $B$. As explained before this requires to work with suitable extensions of the algebras $\Bc$, $\Ec$, $\Ac$. Inspired by Cuntz and Quillen \cite{CQ95}, the basic idea is to replace the algebra of smooth functions $\cinf(B)$ with a Fedosov-type deformation of the algebra of (ordinary) differential forms $\Om(B)$. First we consider the graded space 
\be
\Om(B,E)=\bigoplus_{n\geq 0}\Om^n(B,E)\ ,\qquad \Om^0(B,E)=\cinf(B,E)\ ,
\ee 
of smooth $E$-valued differential forms over $B$. One has an isomorphism of vector spaces $\Om^n(B,E)\cong \cinf(M,\pi^*(\La^nT^*B))$, where $\pi^*(\La^nT^*B)$ is the pullback of the vector bundle $\La^nT^*B\to B$ on the total space of the submersion. $\Om(B,E)$ is a right $\Om(B)$-module, for the usual exterior product of differential forms. The graded algebra of differential forms with values in fiberwise pseudodifferential operators is defined analogously:
\be
\Om(B,CL) = \bigoplus_{n\geq 0}\Om^n(B,CL)\ ,\qquad CL=\bigcup_{k\in\zz} CL^k\ .
\ee
It acts naturally as (left) endomorphisms on the module $\Om(B,E)$. Then we need some extra structure in order to define a connection on $E$. Recall that the set of vertical vector fields on $M$ is the kernel of the tangent map $\pi_*: TM\to TB$, or equivalently the subbundle of $TM$  tangent to the submanifolds $(M/B)_b$, $b\in B$. Then choose a horizontal distribution $H\subset TM$, i.e. a direct summand for the vertical vector fields: $TM=\ker(\pi_*)\oplus H$. This provides a lifting $h:\cinf(B,TB)\to\cinf(M,H)$ of the vector fields from the base to the total space of the submersion. A connection $\nabla: \Om^n(B,E)\to \Om^{n+1}(B,E)$ is given by the usual formula
\beq
(\nabla\xi)(X_0,\ldots ,X_n) &=& \sum_{i=0}^n (-)^i h(X_i)\cdot \xi(X_0,\ldots ,\hat{X_i},\ldots, X_n) \\
&& +\sum_{i<j} (-)^{i+j}\xi([X_i,X_j],X_0,\ldots ,\hat{X_i},\ldots,\hat{X_j},\ldots,X_n)\non
\eeq
for any vector fields $X_0,\ldots ,X_n$ over $B$ and $\xi\in \Om^n(B,E)$. $\nabla$ is a derivation of right $\Om(B)$-module: $\nabla(\xi\om)=(\nabla\xi)\om + (-)^n\xi d\om$ for any $\xi\in \Om^n(B,E)$ and $\om\in\Om(B)$. The curvature $\nabla^2$ is the endomorphism $\theta\in \Om^2(B,CL^1)$ mapping two basic vector fields $X,Y\in \cinf(B,TB)$ to the vertical vector field
\be 
\theta(X,Y) = [h(X),h(Y)] - h([X,Y])\ \in \cinf(M,\ker(\pi_*))\ .
\ee
For any $n\in\nn$ let the quotient $LS^{n} = CL^{n}/CL^{n-1}$ denote the space of leading symbols of order $n$. We introduce the following algebras of operator-valued differential forms with compact support on $B$:
$$
\Rc_0=\bigoplus_{n\geq 0}\Omc^{n}(B,CL^{n-1}) \ ,\ \Mc_0=\bigoplus_{n\geq 0}\Omc^{n}(B,CL^{n})\ ,\ \Pc_0=\bigoplus_{n\geq 0}\Omc^{n}(B,LS^{n})
$$
$\Mc_0$ acts naturally by endomorphisms on $\Om(B,E)$, $\Rc_0$ is a two-sided ideal in $\Mc_0$ and $\Pc_0$ is the quotient algebra. Observe that $\Bc$, $\Ec$, $\Ac$ are exactly the subalgebras of degree zero forms in $\Rc_0$, $\Mc_0$, $\Pc_0$ respectively. Now, the connection acts on endomorphisms by the odd derivation $\delta = [\nabla,\ ]$. It is easy to see that $\delta$ leaves $\Rc_0$ and $\Mc_0$ globally invariant, hence it acts on $\Pc_0$ also. Moreover the curvature $\te\in \Om^2(B,CL^{1})$ is a multiplier of $\Rc_0$ and $\Mc_0$, hence of $\Pc_0$. The derivation $\delta$ is not a differential since $\delta^2=[\te,\ ]\neq 0$. Using a trick of Connes (\cite{C94} pp. 229), we shall enlarge $\Mc_0$ by adding a multiplier $v$ such that $v^2=\te$ and $\om_1 v \om_2=0$ for any $\om_1,\om_2\in\Mc_0$. One can think of $v$ as having form degree one. Hence the resulting graded algebra $\Mc_0[v]$ is the set of elements 
\be
\al = \om_{11} + \om_{12}v + v\om_{21} + v\om_{22}v\ ,\quad \om_{ij}\in \Mc_0\ .
\ee
The crucial fact is that $\Mc_0[v]$ is provided with a differential $d$ of degree one defined by the relations $d\om=\delta\om + v\om +(-)^n\om v$ if $\om\in\Mc_0$ has degree $n$, and $dv=0$. One checks that $d^2=0$ which turns $\Mc_0[v]$ into a differential graded algebra. We denote by $\Mc$ the \emph{even degree} subspace of $\Mc_0[v]$ \emph{endowed with the Fedosov product}
\be
\al_1\odot\al_2 = \al_1\al_2 - d\al_1 d\al_2\ .
\ee
$\Mc$ is an associative (trivially graded) algebra. The map $\Mc\to\cinfc(B,CL^0)$, which projects an element $\al=\om_{11} + \om_{12}v + v\om_{21} + v\om_{22}v$ to its component of degree zero (equivalently the component of degree zero of the differential form $\om_{11}$), defines a linearly split homomorphism $\Mc\to \Ec$. Hence $\Mc$ is an extension of $\Ec$. One proceeds similarly with $\Rc_0$ and $\Pc_0$: the algebras $\Rc$ and $\Pc$ are defined as the even subspaces of $\Rc_0[v]$ and $\Pc_0[v]$ respectively, endowed with the Fedosov product. Again the projections onto the degree zero components induce surjective homomorphisms $\Rc\to\Bc$ and $\Pc\to\Ac$. Moreover the extension $0\to\Rc\to\Mc\to\Pc\to 0$ is a lifting of $(E)$. Putting everything together we have built a diagram of extensions
\be
\vcenter{\xymatrix{
 & 0 \ar[d] & 0 \ar[d] & 0 \ar[d] &  \\
0 \ar [r] & \Jc \ar[r] \ar[d] & \Nc \ar[r] \ar[d] & \Qc  \ar[r] \ar[d] & 0  \\
0 \ar [r] & \Rc \ar[r] \ar[d] & \Mc \ar[r] \ar[d] & \Pc  \ar[r] \ar[d] & 0  \\
0 \ar [r] & \Bc \ar[r] \ar[d] & \Ec \ar[r] \ar[d] & \Ac  \ar[r] \ar[d] & 0  \\
 & 0 & 0 & 0 &  }}
\ee
There is a concrete description of the first line. The ideal $\Nc\subset \Mc$ is the set of elements $\al=\om_{11} + \om_{12}v + v\om_{21} + v\om_{22}v$, with $\om_{11},\om_{22}\in \Mc_0$ of even degree, $\om_{12},\om_{21}\in \Mc_0$ of odd degree, and $\om_{11}$ has no component of degree zero. Hence $\al$ has an overall degree $\geq 2$, which means that the algebra $\Nc$ is \emph{nilpotent}. Similarly with $\Jc$ and $\Qc$. Hence according to Remark \ref{rem} this implies $\Mch=\Mc$ and any trace over $\Rc^{n+1}$, $n\geq 1$, determines a class in $HP^0(\Bc)$. \\

We shall now construct the connecting morphism of the extension $(E)$ as explained in section \ref{scm}. A linear splitting $\si:\Ac\to\Mc$ is obtained as follows: first choose a ``quantization map'' $q:\Ac\to\Ec$, which associates to any function $a\in \cinfc(B,LS^0)$ a family of pseudodifferential operators $q(a)\in \cinfc(B,CL^0)$ with leading symbol $a$. Then map $\Ec$ to the degree zero subspace of $\Mc$. This yields the desired linear splitting $\si$ for the extension
\be
\vcenter{\xymatrix{
0\ar[r] & \Rc + \Nc \ar[r]  & \Mc \ar[r] & \Ac \ar[r] \ar@{.>}@/_1pc/[l]_{\si}  & 0 }}
\ee
whence a classifying homomorphism $\si_*:T\Ac\to\Mc$ which intertwines the tensor product and the Fedosov product: $\si_*(a_1\otimes\ldots\otimes a_n)=\si(a_1)\odot\ldots\odot\si(a_n)$ for any element $a_1\otimes\ldots\otimes a_n\in T\Ac$. By construction $\si_*$ restricts to a homomorphism $J\Ac\to \Rc+\Nc$. To see this, recall that the ideal $J\Ac=\ker(T\Ac\to\Ac)$ is generated by the differences $a_1a_2-a_1\otimes a_2$ for any pair $a_1,a_2\in\Ac$. Then
\beq
\si_*(a_1a_2-a_1\otimes a_2) &=& \si(a_1a_2)-\si(a_1)\odot \si(a_2)\\
&=& \si(a_1a_2)-\si(a_1)\si(a_2) + d\si(a_1)d\si(a_2)\ .\non
\eeq
The first term of the r.h.s. $\si(a_1a_2)-\si(a_1)\si(a_2)\in \cinfc(B,CL^{-1})$ lies in the degree zero subspace of $\Rc$, whereas the two-form $d\si(a_1)d\si(a_2)$ lies in $\Nc$. Hence $J\Ac$ is mapped to the ideal $\Rc+\Nc$ as claimed, and $\si_*$ extends to a homomorphism of pro-algebras (recall that $\Nc$ is nilpotent)
\be
\si_*:\Th\Ac\to\Mct=\varprojlim_n\Mc/(\Rc+\Nc)^n=\varprojlim_n\Mc/\Rc^n
\ee 
Now consider the classes in $HP^0(\Bc)$ represented by traces on the powers of $\Rc$. We shall construct such traces by combining closed currents (cycles) in the base manifold $B$ with the ordinary (fiberwise) trace of pseudodifferential operators. Indeed if $\om\in \Omc(B,CL^k)$ takes values in the space of pseudodifferential operators of order $k< -\dim (M/B)$, the operator trace is well-defined and yields a smooth differential form $\Tr(\om)\in \Omc(B)$. 

\begin{lemma}
Let $C$ be a cycle of dimension $2m$ in $B$. The linear functional $\tau: \Rc^{n+1}\to\cc$, defined for $n=\dim (M/B)+2m$ 
\be
\tau(\om_{11} + \om_{12}v + v\om_{21} + v\om_{22}v) = \frac{m!}{(2m)!}\int_C \Tr(\om_{11} - \om_{22} \te)
\ee
is a trace, i.e. vanishes on the commutator subspace $[\Rc^n,\Rc]$.
\end{lemma}
{\it Proof:} The algebra $\Rc_0$ is the direct sum of the spaces $\Omc^{k}(B,CL^{k-1})$. Let us take the $n$-fold (ordinary) product of operator-valued differential forms:
$$
\Omc^{k_1}(B,CL^{k_1-1}) \times \ldots \times \Omc^{k_n}(B,CL^{k_n-1}) \to \Omc^{k_1+\ldots+k_n}(B,CL^{k_1+\ldots+k_n-n})\ .
$$
Integration over $C$ will retain the case $k_1+\ldots+k_n= \dim C$ only. Hence the pseudodifferential order $k_1+\ldots+k_n-n$ is $\dim C - n$. If moreover one chooses $n>\dim(M/B)+\dim C$, the order is $<-\dim(M/B)$ and the corresponding pseudodifferential operators are trace-class. This estimate does not change if one takes further products with the multiplier $\te\in \Om^2(B,CL^1)$, or if the ordinary product is replaced by the Fedosov product (which increases the form degree by two and the order by at most one). One concludes that the linear functional $\tau$ is well-defined on $\Rc^{n+1}$ provided one chooses $n=\dim(M/B)+\dim C$. \\
Then one checks as in \cite{C94} pp. 229 that $\tau(\al)$ vanishes if $\al$ is the graded commutator of elements in the DG algebra $\Rc_0[v]$, or if $\al=d\beta$ is closed. Hence $\tau$ is a trace for the Fedosov product. \cqfd\\

The factor $m!/(2m)!$ is (up to a sign) the correct normalization needed for passing from the $X$-complex to the de Rham complex (\cite{CQ95}). Thus $\tau$ is a cocycle of even degree over the subcomplex $F^{2n+1}_{\Rc}X(\Mc)$ of the $\Rc$-adic filtration of $X(\Mc)$, provided $n=\dim(M/B)+\dim C$. The next step is to extend $\tau$ to a linear map $\tau_R:\Mc\to\cc$ and view it as a cochain over the whole complex $X(\Mc)$. We use zeta-function renormalization \cite{P6}. Fix a smooth family of fiberwise elliptic positive pseudodifferential operators $D$ of order one. For example, $D$ may be taken as the square root of a fiberwise Laplacian associated to some smooth family of Riemannian metrics on the fibers of the submersion. For $\om\in \Omc(B,CL)$ of any pseudodifferential order, the zeta-function $z\in\cc \mapsto \Tr(\om D^{-z})\in \Omc(B)$ is holomorphic on a half-plane $\re(z)\gg 0$ and admits a meromorphic extension to the entire plane with only simple poles. Taking the \emph{finite part} of this function at $z=0$ thus defines a renormalization of the operator trace:
\be
\Pf_{z=0}\Tr(\om D^{-z}) \in \Omc(B)\ ,\quad \forall \om\in \Omc(B,CL)\ .
\ee

\begin{definition}
Let $C$ be a cycle of dimension $2m$ in $B$ and let $\tau$ be the associated trace over $\Rc^{n+1}$, $n=\dim (M/B)+2m$. Choose an elliptic and positive family of pseudodifferential operators $D\in \cinf(B,CL^1)$. The zeta-function renormalization of $\tau$ is the linear map $\tau_R:\Mc\to\cc$ 
\be
\tau_R(\om_{11} + \om_{12}v + v\om_{21} + v\om_{22}v) = \frac{m!}{(2m)!}\int_C \Pf_{z=0} \Tr\big((\om_{11} - \om_{22} \te)D^{-z}\big)
\ee
 for any $\om_{11} + \om_{12}v + v\om_{21} + v\om_{22}v \in \Mc$.
\end{definition} 

Of course $\tau_R$ is not a trace on $\Mc$ because the insertion of the operator $D^{-z}$ destroys the cyclicity of the operator trace. It is well-known however that the obstruction for the zeta-renormalized trace to be a true trace on the algebra of pseudodifferential operators is expressed in terms of the (fiberwise) Wodzicki residue \cite{Wo}. Indeed if $\om_1,\om_2\in \Omc(B,CL)$ are differential forms with values in pseudodifferential operators of any order, the zeta-renormalized operator trace applied to their graded commutator yields a residue (see e.g. \cite{Pa})
\be
\Pf_{z=0}\Tr([\om_1,\om_2]D^{-z}) = \res \Tr(\om_1[\ln D,\om_2]D^{-z})\ ,
\ee
The logarithm $\ln D$ does not belong to the algebra of classical pseudodifferential operators, but the commutator $[\ln D,\om]$ is a well-defined element of $\Omc(B,CL)$ modulo smoothing operators. Indeed it admits the asymptotic expansion
\be
[\ln D,\om] \sim [D,\om] D^{-1} -\frac{1}{2}[D,[D,\om]] D^{-2} + \frac{1}{3} [D,[D,[D,\om]]] D^{-3} - \ldots
\ee
where the commutator $[D,\ ]$ does not increase the order of operators since $D$ is of order one. The residue $\res\Tr(\om D^{-z})$ is a differential form over $B$ which depends on the complete symbol of $\om$ only and hence kills all smoothing operators: it is the integral, over  the cotangent sphere bundle $S^*M/B$, of the order $-\dim (M/B)$ component in the asymptotic expansion of the symbol \cite{Wo}. We define a linear functional $\Omc(B,CL)\to \cc$ by integration of the fiberwise Wodzicki residue over the cycle $C$:
\be
\bint_C \om := \int_C \res \Tr(\om D^{-z})  = \int_C\int_{S^*M/B} s(\om)_{-n}\  \eta(d\eta)^{n-1}  \label{bint} 
\ee
where $s(\om)$ is the complete symbol of $\om$, $n=\dim (M/B)$, $\eta$ is the canonical one-form on the cotangent bundle $T^*M/B$ of the submersion fibers, and $\int_{S^*M/B}$ denotes integration along the cotangent sphere bundle. It follows from the properties of the Wodzicki residue that (\ref{bint}) does not depend on the choice of (elliptic, positive, order one) operator $D$, and defines a $\delta$-closed graded trace over $\Omc(B,CL)$. This allows to express the boundary of $\tau_R$ viewed as a cochain of even degree over $X(\Mc)$. Indeed the boundary map $\d:\Om^1\Mc_{\nat}\to\Mc$ is given by the Fedosov commutator $\d(\al_1\dd\al_2) = \al_1\odot\al_2 - \al_2\odot\al_1$, so that the composition $\tau_R \d (\al_1\dd\al_2) = \tau_R(\al_1\odot\al_2 - \al_2\odot\al_1)$ has to be a sum of Wodzicki-type residues. For simplicity we shall only evaluate $\tau_R\d$ on the range of the chain map $\si_*: X(T\Ac)\to X(\Mc)$. In odd degree the range is linearly generated by elements of type $\nat((\si_1\odot\ldots\odot\si_n)\dd\si_{n+1})\in \Om^1\Mc_{\nat}$ where $\si_i=\si(a_i)$ for some $a_i\in\Ac$. Writing the Fedosov products by means of differential forms, this is equivalent to the linear span of elements of type $\nat\si_0d\si_1\ldots d\si_{2n}\dd\si_{2n+1}$ and $\nat d\si_1\ldots d\si_{2n}\dd\si_{2n+1}$. Notice that the derivative $\delta\ln D=[\nabla,\ln D]$ is a pseudodifferential operator with asymptotic expansion
\be
\delta\ln D \sim \delta D D^{-1} -\frac{1}{2}[D,\delta D] D^{-2} + \frac{1}{3} [D,[D,\delta D]] D^{-3} - \ldots
\ee

\begin{proposition}\label{pfor}
Let $\nat\si_0d\si_1\ldots d\si_{2n}\dd\si_{2n+1}$ and $\nat d\si_1\ldots d\si_{2n}\dd\si_{2n+1}$ be generic odd chains in $X(\Mc)$ such that all $\si_i$'s are in the image of the linear splitting $\si:\Ac\to\Mc$. Let $C$ be a cycle of even dimension with associated trace $\tau$. Then the boundary of the renormalized trace $\tau_R$ reads
\beq
\lefteqn{\tau_R\d(\si_0d\si_1\ldots d\si_{2n}\dd\si_{2n+1}) =  \frac{n!}{(2n)!}\bint_C (\si_0d\si_1\ldots d\si_{2n})_{11}[\ln D,\si_{2n+1}] } \label{form1}\\
&&\qquad \qquad + \frac{(n+1)!}{(2n+2)!}\bint_C  (\si_0d\si_1\ldots d\si_{2n+1} - \si_{2n+1}d\si_0\ldots d\si_{2n})_{11}\delta\ln D \non
\eeq
\be
\tau_R\d(d\si_1\ldots d\si_{2n}\dd\si_{2n+1}) =  \frac{n!}{(2n)!}\bint_C (d\si_1\ldots d\si_{2n})_{11}[\ln D,\si_{2n+1}] \label{form2}
\ee
where for any $\al=\om_{11} + \om_{12}v + v\om_{21} + v\om_{22}v\in M_0[v]$ the bracket $(\al)_{11}$ means projection onto the component $\om_{11}$. If $k=2n$ one finds
\beq
(\si_0d\si_1\ldots d\si_k)_{11} &=& \si_0\delta\si_1\ldots\ldots\ldots \delta\si_k\non\\
&+& \sum_{i=1}^{k-1} \si_0\delta\si_1\ldots \si_i\te\si_{i+1} \ldots \delta\si_k \label{toto}\\
&+&  \sum_{i=1}^{k-3}\sum_{j=i+2}^{k-1} \si_0\delta\si_1\ldots \si_i\te\si_{i+1}\ldots \si_j\te\si_{j+1} \ldots \delta\si_k\non\\
&&  \vdots \non\\
&+& \si_0 \si_1\te\si_2\ldots \si_{k-1}\te\si_k \ ,\non
\eeq
whereas if $k=2n+1$ the last line is $\sum_{i=1}^{k}\si_0 \si_1\te\si_2\ldots\delta\si_i\ldots \si_{k-1}\te\si_k$. Similarly with $(d\si_1\ldots d\si_{2n})_{11}$.
\end{proposition}
{\it Proof:} By definition the boundary map $\d=\bb:\Om^1\Mc_{\nat}\to\Mc$ carries an element $\nat \al\dd \beta$ to the Fedosov commutator $\al\odot\beta - \beta\odot\al= [\al,\beta]-d\al d\beta + d\beta d\al$. Therefore
\beq
\tau_R\d(\si_0d\si_1\ldots d\si_{2n}\dd\si_{2n+1}) &=& \tau_R([\si_0d\si_1\ldots d\si_{2n},\si_{2n+1}]) \non\\
&& + \tau_R(d\si_{2n+1}d\si_0\ldots d\si_{2n}-d\si_0\ldots d\si_{2n+1}) \non\\
\tau_R\d(d\si_1\ldots d\si_{2n}\dd\si_{2n+1}) &=& \tau_R([d\si_1\ldots d\si_{2n},\si_{2n+1}])\ .\non
\eeq
Let us first evaluate $\tau_R$ on a commutator $[\al,\si]$ where $\al=\om_{11} + \om_{12}v + v\om_{21} + v\om_{22}v$ is an element of degree $2n$ and $\si$ is of form degree zero. One has $\al\si=\om_{11}\si + v\om_{21}\si$ and $\si\al=\si\om_{11}+ \si\om_{12}v$, hence
$$
\tau_R([\al,\si])=\frac{n!}{(2n)!}\int_C \Pf_{z=0}\Tr([\om_{11},\si]D^{-z})=\frac{n!}{(2n)!} \bint_C \om_{11}[\ln D,\si]\ .
$$
Applying this to the forms $\al=\si_0d\si_1\ldots d\si_{2n}$ or $\al=d\si_1\ldots d\si_{2n}$ and $\si=\si_{2n+1}$ yields the first terms in (\ref{form1}, \ref{form2}). Then we evaluate $\tau_R$ on a coboundary $d\al$ where $\al=\om_{11}+\om_{12}v$ is an odd element of form degree $2n+1$. One has $d\al=\delta\om_{11}+\om_{12}\te +(\delta\om_{12}-\om_{11})v+v\om_{11}+v\om_{12}v$, so that 
\beq
\tau_R(d\al) &=& \frac{(n+1)!}{(2n+2)!}\int_C \Pf_{z=0}\Tr\big((\delta\om_{11}+\om_{12}\te-\om_{12}\te)D^{-z}\big)\non\\
&=& \frac{(n+1)!}{(2n+2)!}\int_C \Pf_{z=0}\Tr(\om_{11}\delta D^{-z})\non\\
&=& -\ \frac{(n+1)!}{(2n+2)!}\bint_C \om_{11} \delta\ln D \non
\eeq
where we used an integration by parts in the second equality (remark that the form $\om_{11}$ is odd), and the third equality can be found for example in \cite{Pa}. Applying this to the form $\al=\si_{2n+1}d\si_0\ldots d\si_{2n}-\si_0d\si_1\ldots d\si_{2n+1}$ yields the second term in (\ref{form1}). Formula (\ref{toto}) is straightforward using $d\si_i=\delta\si_i+v\si_i+\si_iv$ and $v^2=\te$.\cqfd\\

One knows that $\tau_R\d$ vanishes on the subcomplex $F^{n+1}X(\Mc)$ for $n =\dim (M/B)+\dim C$, hence extends to a cocycle over $X(\Mct)$. Finally the cyclic cohomology class $E^*([\tau])\in HP^1(\Ac)$ is represented by the composition of chain maps 
\be
X(\Th\Ac)\stackrel{\si_*}{\longrightarrow} X(\Mct) \stackrel{\tau_R\d}{\longrightarrow} \cc\ . 
\ee
Here we can interpret the fact that $\tau_R\d$ extends to a cocycle over $X(\Mct)$ by the property of the fiberwise Wodzicki residue that ignores the pseudodifferential operators of low order, that is, the high powers of the ideal $\Rc\subset \Mc$. It remains to evaluate the pairing of $E^*([\tau])$ with an elliptic symbol class $[g]\in K_1(\Ac)$. Let $g\in GL_{\infty}(\Ac)$ be a representative of $[g]$. For notational simplicity we shall forget the stabilization by matrices and suppose that $g\in GL_1(\Ac)\subset \Ac^+$. Then $\si(g)$ and $\si(g^{-1})$ are two families of elliptic pseudodifferential operators over $B$ such that 
\be
1-\si(g)\si(g^{-1})\in \cinfc(B,CL^{-1})\subset \Rc\ .\label{comp}
\ee
Let $Q$ be the image of $\si(g)$ under the natural map $\Mc^+\to \Mct^+$. Then $Q\in GL_1(\Mct)$. Indeed, $1-\si(g)\odot\si(g^{-1})=1-\si(g)\si(g^{-1}) +d\si(g)d\si(g^{-1})$ is in the ideal $\Rc+\Nc$ by virtue of (\ref{comp}), and one can easily show that the inverse of $Q$ (for the Fedosov product) is given by the series
\be
Q^{-1}=\sum_{n=0}^{\infty} \si(g^{-1})\odot(1-\si(g)\odot\si(g^{-1}))^{\odot n}\ \in \Mct^+\ .
\ee
There is an equivalent description of the Fedosov inverse involving the parametrix $P\equiv Q^{-1}\mod \Nc$ of $Q$. This allows to write $Q^{-1}$ in terms of differential forms:
\be
Q^{-1} = \sum_{n=0}^{[\frac{\dim B}{2}]} P(dQdP)^n\ ,\quad P = \sum_{n=0}^{\infty} \si(g^{-1})(1-\si(g)\si(g^{-1}))^n\ .
\ee
Now we can calculate the pairing as  $\langle E^*([\tau]),[g] \rangle = \frac{1}{\sqrt{2\pi i}}\tau_R\d(Q^{-1}\dd Q)$. Combining this with the explicit formula of Proposition \ref{pfor} and the Index Theorem \ref{t} gives an expression of $\langle [\tau], \Ind_E([g]) \rangle$ in terms of $Q$, its parametrix $P$ and the fiberwise Wodzicki residue. The computation is tedious but straightforward. We shall state the result in an elegant way using Chern-Simons forms. Introduce an infinitesimal parameter $\eps$ of odd degree, which means $\eps^2=0$. The superconnection
\be
\nabla^{\eps}_D := \nabla +\eps \ln D
\ee
acts on the algebra $\Omc(B,CL)[\eps]=\Omc(B,CL)\oplus \eps\, \Omc(B,CL)$ by graded commutator. Its curvature is $(\nabla +\eps \ln D)^2 = \te -\eps\delta\ln D$. The ``adjoint'' action of $Q$ gives a new superconnection $P \nabla^{\eps}_D Q$. Now if $\nabla_0$ and $\nabla_1$ are two superconnections, we let $\nabla_t=(1-t)\nabla_0+ t\nabla_1$ be the linear interpolation for $t\in [0,1]$. The associated Chern-Simons form is the following element of even degree in $\Omc(B,CL)$ (defined as always modulo smoothing operators, due to the presence of $\ln D$):
\be
\cs (\nabla_0,\nabla_1) = \int_0^1 dt\, (\nabla_1-\nabla_0) e^{\nabla_t^2}|_{\eps}\ ,
\ee
where $|_{\eps}$ means that we only take the $\eps$-component in $\Omc(B,CL)[\eps]$. Since $\nabla$ is of form degree one and $\eps$ is nilpotent, the exponential is actually a polynomial in the curvature $\nabla_t^2$. Applying this to $\nabla_0=\nabla^{\eps}_D$ and $\nabla_1=P \nabla^{\eps}_D Q$ one has

\begin{corollary}
Let $M\to B$ be a proper submersion with connection $\nabla$. Let $[g]\in K_1(\Ac)$ be the symbol class of an elliptic family of fiberwise pseudodifferential operators $Q$ with parametrix $P$. Let $C$ be a cycle of even dimension in the base manifold $B$, and $[\tau] \in HP^0(\Bc)$ the associated cyclic cohomology class. The evaluation of $[\tau]$ on the index class $\Ind_E([g])\in K_0(\Bc)$ is given by the fiberwise residue 
\be
\langle [\tau], \Ind_E([g]) \rangle = \bint_C \cs (\nabla^{\eps}_D,P\nabla^{\eps}_DQ)
\ee
where $\nabla_D^{\eps}$ is the superconnection $\nabla + \eps\ln D$, and $D$ is any family of elliptic positive pseudodifferential operators of order one.  \cqfd\\
\end{corollary}

Let us display some useful formulas in low dimension. If $C$ is just a point in the base manifold $B$, the above pairing calculates the index of the elliptic operator $Q$ at point $C$. The formula amounts to a Wodzicki residue on one fiber:
\be
\langle [\tau], \Ind_E([g]) \rangle = \bint_C P[\ln D,Q]\ . \label{rad}
\ee
One recognizes the Radul cocycle \cite{Ra} evaluated on $Q$ and its parametrix. It is instructive to check that the number (\ref{rad}) does not depend on the choice of $D$. Indeed if $D'$ is another elliptic operator of order one, the difference $\ln D'-\ln D$ is a classical pseudodifferential operator, hence the commutator with $\ln D'-\ln D$ is an inner derivation. It follows from the trace property of the Wodzicki residue that (\ref{rad}) remains unchanged. \\
Now if $C$ is a two-dimensional cycle then the above pairing computes the evaluation of $C$ on the first Chern class of the index bundle associated to the elliptic family $Q$. One obtains
\beq
\lefteqn{\langle [\tau], \Ind_E([g]) \rangle = }\\
&& \frac{1}{2} \bint_C\big( P\delta Q\delta P[\ln D,Q] + P(\delta Q \delta\ln D -\delta\ln D\delta Q)  + (\te P+P\te)[\ln D,Q] \big) \non
\eeq
In higher dimensions the formulas involve increasing powers of $\delta Q\delta P$ and of the curvature $\te$.


\begin{thebibliography}{9}

\bibitem{C86} A. Connes: Non-commutative differential geometry, {\it Publ. Math. IHES} {\bf 62} (1986) 41-144.

\bibitem{C94} A. Connes: {\it Non-commutative geometry}, Academic Press, New-York (1994).

\bibitem{CM95} A. Connes, H. Moscovici: The local index formula in non-commutative geometry, GAFA {\bf 5} (1995) 174-243.

\bibitem{Cu89} J. Cuntz: Universal extensions and cyclic cohomology, {\it C. R. Acad. Sci. Paris} {\bf 309} S{\'e}rie I (1989) 5-8.

\bibitem{CQ95} J. Cuntz, D. Quillen: Cyclic homology and nonsingularity, JAMS {\bf 8} (1995) 373-442.

\bibitem{CQ97} J. Cuntz, D. Quillen: Excision in bivariant periodic cyclic cohomology, {\it Invent. Math.} {\bf 127} (1997) 67-98.

\bibitem{M} J. Milnor: {\it Algebraic K-theory}, Ann. Math. Studies {\bf 72} Princeton University press (1974).

\bibitem{Ni} V. Nistor: Higher index theorems and the boundary map in cyclic cohomology, {\it Doc. Math. J. DMV} {\bf 2} (1997) 263-295.

\bibitem{Pa} S. Paycha, S. Scott: Chern-Weil forms associated with superconnections, {\it Analysis, geometry and topology of elliptic operators}, World Sci. Publ., Hackensack, NJ (2006) 79-104

\bibitem{P4} D. Perrot: Anomalies and noncommutative index theory, lectures given at Villa de Leyva, Colombia (2005), S. Paycha and B. Uribe Ed., {\it Contemp. Math.} {\bf 434} (2007) 125-160.

\bibitem{P5} D. Perrot: Secondary invariants for Fr\'echet algebras and quasihomomorphisms, {\it Doc. Math. J. DMV} {\bf 13} (2008) 275-363.

\bibitem{P6} D. Perrot: Quasihomomorphisms and the residue Chern character, preprint arXiv:0804.1048.

\bibitem{P7} D. Perrot: Localization over complex-analytic groupoids and conformal renormalization, {\it J. Noncommut. Geom.} {\bf 3} (2009) 289-325.

\bibitem{Ra} A. O. Radul: Lie algebras of differential operators, their central extensions and $W$-algebras.  (Russian) {\it Funktsional. Anal. i Prilozhen.} {\bf 25} (1991) 33-49;  translation in {\it Funct. Anal. Appl.} {\bf 25} (1991) 25-39 

\bibitem{Wo} M. Wodzicki: {\it Non-commutative residue}, Lect. Notes in Math. {\bf 1283}, Springer Verlag (1987).






\end{thebibliography}
\end{document}